\documentclass[11pt,a4paper]{amsart}
\usepackage[dvips]{epsfig}
\usepackage{a4wide}
\usepackage{amscd}
\usepackage{amssymb}
\usepackage{amsthm}
\usepackage{amsmath}
\usepackage{latexsym}
\usepackage{enumerate}
\usepackage{color}

{%
\setcounter{enumi}{0}

\begin{enumerate}}%
{\end{enumerate}
}
{%
\setcounter{enumi}{0}

\begin{enumerate}}%
{\end{enumerate}
}

\newtheorem{lemma}{Lemma}
\newtheorem{remark}{Remark}
\newtheorem{theorem}{Theorem}
\newtheorem{prop}{Proposition}
\newtheorem{cor}{Corollary}
\newtheorem{definition}{Definition}

\newcommand{\charu}{{1\hspace*{-3.5pt}1}} 
\newcommand{\R}{\mathbb{R}}               
\newcommand{\G}{\textbf{g}}               
               
\newcommand{\N}{\mathbb{N}}               
\newcommand{\h}{\delta t}                 
                 
\newcommand{\K}{K}                        
\newcommand{\abs}[1]{\ensuremath{\left|#1\right|}}
\newcommand{\Grad}{\mathrm{\nabla}} 
\newcommand{\disp}{\ensuremath{\displaystyle}}
\newcommand{\dm}{\ell}
\newcommand{\V}{\textbf{V}}              
\newcommand{\ga}{\gamma}                  
\newcommand{\pc}{\mathcal{P}_c}        
         
\newcommand{\E}{\mathcal{E}}
\newcommand{\D}{\mathcal{D}}
\newcommand{\T}{\mathcal{T}}
\newcommand{\B}{\mathcal{B}}
\newcommand{\norm}[1]{\ensuremath{\left\|#1\right\|}}
\newcommand{\di}{\mathrm{div}}
\newcommand{\dd}{\mathrm{d}}
\newcommand{\1}{g} 
\newcommand{\3}{l} 
\newcommand{\2}{l} 
 
\newcommand{\sizt}{h} 
\newcommand{\bord}{\Gamma} 
\newcommand{\pas}{\delta t}

\newcommand{\KL}{K|L}
\newcommand{\de}{\delta_{K|L}^{n+1}}
\newcommand{\dis}{d_{K|L}}
\newcommand{\sig}{\sigma_{K|L}}
\newcommand{\tokl}{\tau_{K|L}}
\newcommand{\QKn}{Q_K^n}

\newcommand{\sn}{\sum_{n=0}^{N-1}}

\newcommand{\snk}{\sum_{n=0}^{N-1}\sum_{K\in \mathcal{T}}}
\newcommand{\sm}{\sn\delta t\sum_{K\in \mathcal{T}}\sum_{L\in N(K)}}
\newcommand{\Lll}{L}

\newcommand{\hh}{\mathcal{D}_m} 
\newcommand{\hhm}{(\mathcal{D}_m)_m} 
\newcommand{\m}{{m\in \N}} 
\newcommand{\ch}{_{\mathcal{D}_m}}
\newcommand{\eh}{_{\mathcal{D}_m}}

\newcommand{\slmin}{s_{l,min}^{n+1}}
\newcommand{\size}{\text{size}}

\newcommand{\sli}{(s^I_{l,K})^{n+1}}
\newcommand{\sgi}{(s^I_{g,K})^{n+1}}

\numberwithin{equation}{section}

\allowdisplaybreaks

\title[Finite volume scheme for compressible]
{Study of  full implicit petroleum engineering finite volume scheme for compressible 
    two phase flow in porous media}

\author[B. Saad]{Bilal Saad}
\address[Bilal Saad]{\newline
         Université de Nantes\newline
         Laboratoire de Math\'ematiques Jean Leray (UMR 6629 CNRS)\newline
         F-44322 Nantes Cedex 032, France   }
               
\email[]{bilal.saad@univ-nantes.fr}

\author[M. Saad]{Mazen Saad}
\address[Mazen Saad]{\newline
         Ecole Centrale de Nantes\newline
         D\'epartement d' Informatique et Math\'ematiques\newline
         Laboratoire de Math\'ematiques Jean Leray (UMR 6629 CNRS)\newline
         1, rue de la No\'e, BP 92101, France   }
         
\email[]{Mazen.Saad@ec-nantes.fr}


\thanks{} 

\begin{document}


\begin{abstract}
  An industrial scheme, to simulate the two compressible phase flow in
  porous media, consists in a finite volume method together with a
  phase-by-phase upstream scheme.  The implicit finite volume scheme
  satisfies industrial constraints of robustness.  We show that the
  proposed scheme satisfy the maximum principle for the saturation, a
  discrete energy estimate on the pressures and a function of the
  saturation that denote capillary terms. These stabilities results
  allow us to derive the convergence of a subsequence to a weak
  solution of the continuous equations as the size of the
  discretization tends to zero. The proof is given for the complete
  system when the density of the each phase depends on the own
  pressure.
\end{abstract}

%
\maketitle
\section{Introduction}\label{sec:introduction}
A rigorous mathematical study of a petroleum engineering schemes takes
an important place in oil recovery engineering for production of
hydrocarbons from petroleum reservoirs. This important problem renews
the mathematical interest in the equations describing the multi-phase
flows through porous media. The derivation of the mathematical
equations describing this phenomenon may be found in \cite{Bear67},
\cite{chavent}. The differential equations describing the flow of two
incompressible, immiscible fluids in porous media have been studied in
the past decades. Existence of weak solutions to these equations has
been shown under various assumptions on physical data \cite{arbogast,
  chavent, chen99, chen2002, chen, feng, GM96, kroener84, kruzkov77}.

The numerical discretization of the two-phase incompressible
immiscible flows has been the object of several studies, the
description of the numerical treatment by finite difference scheme may be
found in the books \cite{aziz}, \cite{peaceman}. 

The finite volume methods have been proved to be well adapted to
discretize conservative equations and  have been
used in industry because they are cheap, simple to code and robust.
The porous media problems are one of the privileged field of
applications.  This success induced us to study and prove the
mathematical convergence of a classical finite volume method for a
model of two-phase flow in porous media. 

For the two-phase incompressible immiscible flows, the convergence of
a cell-centered finite volume scheme to a weak solution is studied in
\cite{michel2003}, and for a cell-centered finite volume scheme, using
a ``phase by phase'' upstream choice for computations of the fluxes
have been treated in \cite{Eymard00} and in \cite{brenier}. The authors
give an iterative method to calculate explicitly the phase by phase
upwind scheme in the case where the flow is driven by gravitational
forces and the capillary pressure is neglected. An introduction of the
cell-centered finite volume can be found in \cite{Eymard:book}.

For the convergence analysis of an approximation to miscible fluid
flows in porous media by combining mixed finite element and finite
volume methods, we refer to \cite{amaziane-ossmani}, \cite{amirat-fv}.


Pioneers works have been done recently by C. Galusinski and M. Saad in
a serie of articles about ``Degenerate parabolic system for
compressible, immiscible, two-phase flows in porous media''
(\cite{CS04}, \cite{CS08}, \cite{CS07}) when the densities depend on
the global pressure , and by Z. Khalil and M. Saad in (\cite{ZS10},
\cite{ZS11}) for the general case where the density of each phase
depends on its own pressure. And for the two compressible, partially
miscible flow in porous media, we refer to \cite{CS10},
\cite{these-bilal}. For the convergence analysis of a finite volume
scheme for a degenerate compressible and immiscible flow in porous
media with the feature of global pressure, we refer to \cite{ZS11fv}.

In this paper, we consider a two-phase flow model where the fluids are
immiscible. The medium is saturated by a two compressible phase flows.
The model is treated without simplified assumptions on the density of
each phase, we consider that the density of each phase depends on its
corresponding pressure. It is well known that equations arising from
multiphase flow in porous media are degenerated. The first type of
degeneracy derives from the behavior of relative permeability of each
phase which vanishes when his saturation goes to zero. The second type
of degeneracy is due to the time derivative term when the saturation
of each phase vanishes.

This paper deals with construction and convergence analysis of a
finite volume scheme for two compressible and immiscible flow in
porous media without simplified assumptions on the state law of the
density of each phase.

The goal of this paper is to show that the approximate solution
obtained with the proposed upwind finite volume scheme
\eqref{sys_disc:pl}--\eqref{sys_disc:pg} converges as the mesh size
tends to zero, to a solution of system \eqref{eq:conservation_masse}
in an appropriate sense defined in section
\ref{sec:mathematical_formulation}. In section \ref{sec:FV}, we
introduce some notations for the finite volume method and we present
our numerical scheme and the main theorem of convergence.

In section \ref{sec:fundamental_lemma}, we derive three preliminary fundamental lemmas. In fact, we will see that we can't control the discrete gradient of pressure since the mobility of
each phase vanishes in the region where the phase is missing. 
So we are going to use the feature of global pressure. 
 We show that the control of velocities ensures the control of the global
pressure and a dissipative term on saturation in the whole domain regardless of the
presence or the disappearance of the phases. \\
 Section \ref{sec:basic-apriori} is devoted to a maximum principle on saturation and  a well posedness of the scheme which inspired from  H.W. Alt, S. Luckhaus \cite{Alt83}. Section
\ref{sec:compacity} is devoted to a space-time $L^1$ compactness of
sequences of approximate solutions. \\
Finally, the
passage to the limit on the scheme and convergence analysis are  performed in section
\ref{sec:limite}. Some numerical results are stated in the last section \ref{sec:numerical}.

\section{Mathematical formulation of the continuous problem}\label{sec:mathematical_formulation}

Let us state the physical model describing the immiscible displacement
of two compressible fluids in porous media. Let $T>0$
be the final time fixed, and let be $\Omega$ a bounded open subset of
$\R^\dm\ (\dm\geq1)$. We set $Q_T=(0,T)\times \Omega$, $\Sigma_T =
(0,T)\times \partial \Omega$. The mass conservation of each phase is given in $Q_T$ 
\begin{equation}\label{eq:conservation_masse}
  \phi(x)\partial_{t}( \rho_{\alpha}(p_\alpha)s_\alpha)(t,x) + \di (\rho_{\alpha}(p_\alpha) \V_{\alpha})(t,x)
  +\rho_\alpha(p_\alpha) s_\alpha  f_{P}^{~}(t,x) = \rho_\alpha(p_\alpha) s^I_\alpha f_{I}^{~}(t,x),
\end{equation}
where $\phi$, $\rho_\alpha$ and $s_\alpha$ are respectively the
porosity of the medium, the density of the $\alpha$ phase and the
saturation of the $\alpha$ phase. Here the functions $f_{I}^{~}$ and
$f_{P}^{~}$ are respectively the injection and production terms. Note
that in equation \eqref{eq:conservation_masse} the injection term is
multiplied by a known saturation $s^I_\alpha$ corresponding to the
known injected fluid, whereas the production term is multiplied by the
unknown saturation $s_\alpha$
corresponding to the produced fluid.\\
The velocity of each fluid $\V_\alpha$ is given by the Darcy law:
\begin{equation}
  \V_{\alpha}= - {\bf K}
  \frac{k_{r_\alpha}(s_\alpha)}{\mu_{\alpha}}\big( \nabla
  p_{\alpha}-\rho_\alpha(p_\alpha){\bf g}\big),\qquad\quad \alpha
  = \2, \1.
\end{equation}
where ${\bf K}$ is the permeability tensor of the porous medium, $k_{r_\alpha}$ the
relative permeability of the $\alpha$ phase, $\mu_\alpha$ the constant
$\alpha$-phase's viscosity, $p_\alpha$ the $\alpha$-phase's pressure
and ${\bf g }$ is the gravity term.
Assuming that the phases occupy the whole pore space, the phase
saturations satisfy
\begin{equation}\label{def:saturation}
  s_{\2}+ s_{\1} = 1.
\end{equation}
The curvature of the contact surface between the two fluids links the
jump of pressure of the two phases to the saturation by the capillary
pressure law in order to close the system
\eqref{eq:conservation_masse}-\eqref{def:saturation}
\begin{equation}\label{def:pression_capillaire.}
  p_c(s_\3(t,x)) = p_{\1}(t,x) - p_{\2}(t,x).
\end{equation}
With the arbitrary choice of \eqref{def:pression_capillaire.} (the jump
of pressure is a function of $s_\3$), the application $s_\3\mapsto
p_c(s_\3)$ is non-increasing, $(\frac{\dd p_c}{\dd s_\3}(s_\3) < 0,
\mbox{ for all } s_\3 \in [0,1])$, and usually $p_c(s_\3=1)=0$
 when the wetting fluid is at its maximum
saturation. 

\subsection{Assumptions and main result}
The model is treated without simplified assumptions on the density of
each phase, we consider that the density of each phase depends on its
corresponding pressure. The main point is to handle a priori estimates
on the approximate solution. The studied system represents two kinds
of degeneracy: the degeneracy for evolution terms
$\partial_t(\rho_\alpha s_\alpha)$ and the degeneracy for dissipative
terms $\di(\rho_\alpha M_\alpha\nabla p_\alpha)$ when the saturation
vanishes. We will see in the section \ref{sec:basic-apriori} that we
can't control the discrete gradient of pressure since the mobility of
each phase vanishes in the region where the phase is missing. So, we
are going to use the feature of global pressure to obtain uniform
estimates on the discrete gradient of the global pressure and the
discrete gradient of the capillary term ${\mathcal B}$ (defined on
\eqref{def:beta}) to treat the degeneracy of this system.

Let us summarize some useful notations in the sequel. We recall the
conception of the global pressure as describe in \cite{chavent}
$$
M(s_\3)\nabla p = M_\2(s_\2) \nabla p_\2 + M_\1(s_\1) \nabla p_\1,
$$
with the $\alpha$-phase's mobility $M_\alpha$ and the total mobility are defined by 
$$
M_{\alpha}(s_{\alpha})=k_{r_\alpha}(s_{\alpha})/ \mu_{\alpha}, \quad
M(s_\3) = M_\2(s_\2)+M_\1(s_\1).
$$
This global pressure $p$ can be written as
\begin{align}\label{def:pression_globale}
  p=p_\1+\tilde{p}(s_\3)=p_\2+\bar{p}(s_\3),
\end{align}
or the artificial pressures are denoted by $\bar{p}$ and $\tilde{p}$
defined by:
\begin{align}\label{def:terme_capillaire}
  \tilde{p}(s_\3)=-\int_{0}^{s_\3} \frac{M_{\2}(z)}{M(z)}
  p_c^{'}(z)\dd z \text{ and } \overline{p}(s_\3)=\int_{0}^{s_\3}
  \frac{M_{\1}(z)}{M(z)} p_c^{'}(z)\dd z
\end{align}
We also define the capillary terms by
$$
\ga (s_\3)=-\frac{M_{\2}(s_\2)M_{\1}(s_\1)}{M(s_\3)}\frac{\dd p_c}{\dd
  s_\3}(s_\3)\geq 0,
$$
and let us finally define the function $\B$ from $[0,1]$ to $\R$
by:
\begin{align}\label{def:beta}
  \mathcal{B}(s_\3)&=\int_{0}^{s_\3}\ga(z) \dd z  =-
  \int_{0}^{s_\3}\frac{M_\2(z)M_\1(z)}{M(z)}\frac{\dd p_c}{\dd
    s_\3}(z) \dd z \notag\\
     & =- \int_{0}^{s_\1} M_\2(z)\frac{\dd
    \bar{p}}{\dd s_\3}(z) \dd z  = \int_{0}^{s_\3}
  M_\1(z)\frac{\dd \tilde{p}}{\dd s_\3}(z) \dd z.
\end{align}
We complete the description of the model \eqref{eq:conservation_masse}
by introducing boundary conditions and initial conditions.  To the system
\eqref{eq:conservation_masse}--\eqref{def:pression_capillaire.} we add
the following mixed boundary conditions. We consider the boundary
$\partial \Omega=\Gamma_\3\cup \Gamma_{\emph{imp}}$, where $\Gamma_\2$
denotes the water injection boundary and $\Gamma_{\emph{imp}}$ the
impervious one.
\begin{equation}\label{cd:bord.}
  \left\{ 
    \begin{aligned}      
      p_\2(t,x) = p_\1(t,x)=0 & \text{ on } (0,T)\times\Gamma_\2,  \\
      \rho_\2 \V_\2 \cdot \textbf{n} = \rho_\1 \V_\1 \cdot \textbf{n}
      = 0 & \text{ on } (0,T)\times\Gamma_{imp},
    \end{aligned} 
  \right.
\end{equation}
where $\textbf{n}$ is the outward normal to $\Gamma_{imp}$.

The initial conditions are defined on pressures
\begin{equation}\label{cd:initial.}
  p_{\alpha}(t=0) =p^{0}_{\alpha} \text{ for } \alpha=\2,\1  \text{ in } \Omega
\end{equation}
We are going to construct a finite volume scheme on orthogonal
admissible mesh, we treat here the case where $$K=k\mathcal{I}_d$$
where $k$ is a constant positive. For clarity, we take $k=1$ which
equivalent to change the scale in time.\\
Next we introduce some physically relevant assumptions on the
coefficients of the system.
\begin{enumerate}[({${H}$}1)]
\item \label{hyp:H1} There is two positive constants $\phi_{0}$ and $\phi_{1}$ such that $\phi_{0}\leq
  \phi(x)\leq \phi_{1}$ almost everywhere $x\in \Omega$.
\item \label{hyp:H2} The functions $M_\2$ and $M_\1$ belongs to ${\mathcal
    C}^{0}([0,1],\R^{+})$, $ M_{\alpha}(s_{\alpha}=0)=0.$ In addition,
  there is a positive constant $m_{0}>0$ such that for all $s_\3\in
  [0,1]$,
  $$
  M_\2(s_\2) + M_\1(s_\1)\geq m_{0}.
  $$
\item \label{hyp:H3} $(f_{P}^{~},f_{I}^{~})\in (L^2(Q_T))^2$,
  $f_{P}^{~}(t,x)$, $f_{I}^{~}(t,x) \ge 0$
  almost everywhere $(t,x)\in Q_T$.
\item \label{hyp:H4} The density $\rho_{\alpha}$ is ${\mathcal
    C}^{1}(\R)$, increasing and there exist two positive constants
  $\rho_{m}>0$ and $\rho_{M}>0 $ such that $0<\rho_{m}\leq
  \rho_{\alpha}(p_{\1})\leq \rho_{M}.$
\item \label{hyp:H5} The capillary pressure fonction $p_c\in
  \mathcal{C}^{1}([0,1];\R^{+})$, decreasing and there exists
  $\underline{p_c}>0$ such that $0<\underline{p_c}\leq |\frac{\dd
    p_c}{\dd s_\3}|$.
\item \label{hyp:H6} The function $\ga \in C^{1}\left([0,1];\R^{+}
  \right)$ satisfies $\ga(s_\3)>0$ for $0<s_\3< 1$ and
  $\ga(s_\3=1)=\ga(s_\3=0)=0.$ We assume that $\mathcal{B}^{-1}$ (the
  inverse of $\mathcal{B}(s_\3)=\int_{0}^{s_\3}\ \ga(z) \dd z$) is an
  H\"{o}lder\footnote{This means that there exists a positive constant
    $c$ such that for all $a, b \in [0,\mathcal{B}(1)],$ one has
    $|\mathcal{B}^{-1}(a)-\mathcal{B}^{-1}(b)|\leq c|a -
    b|^{\theta}$.}  function of order $\theta$, with $0<\theta\leq 1,
  \text{ on } [0,\mathcal{B}(1)]$.
\end{enumerate} 
The assumptions ({H}\ref{hyp:H1})--({H}\ref{hyp:H6}) are classical for
porous media. Note that, due to the boundedness of the capillary
pressure function, the functions $\tilde{p}$ and $\bar{p}$ defined in
\eqref{def:terme_capillaire} are bounded on $[0,1]$.

Let us define the following Sobolev space
$$
H^{1}_{\Gamma_\2}(\Omega) = \{ u\in H^{1}(\Omega); u=0 \ sur \
\Gamma_\2 \},
$$
this is an Hilbert space with the norm $\|u
\|_{H^{1}_{\Gamma_\2}(\Omega)} = \|\nabla u\|_{(L^{2 }(\Omega))^{\dm}}$.

\begin{definition}\label{def}$\left(\text{Weak solutions} \right).$
  Under assumptions ({H}\ref{hyp:H1})-({H}\ref{hyp:H6}) and
  definitions \eqref{def:pression_globale}-\eqref{cd:initial.} with
  the fact that $ p^{0}_\2,\ p^{0}_\1 $ belongs to $L^2(\Omega)$ and
  $s^0_\alpha$ satisfies $0\le s^0_\alpha\le 1$ almost everywhere in
  $\Omega$, then the pair $\left(p_\2, p_\1 \right)$ is a weak
  solution of problem \eqref{eq:conservation_masse} satisfying :
\begin{align}
  & p_\alpha \in L^2(0,T;L^2(\Omega)),~
  \sqrt{M_\alpha(s_\alpha)}\nabla p_\alpha \in
  (L^2(0,T;L^2(\Omega)))^\dm, \\
  & 0\leq s_\alpha(t,x)\leq 1 \text{ a.e in } Q_T \; (\alpha=\2,\1),~
  \mathcal{B}(s_\2)\in L^{2}(0,T;H^{1}_{\Gamma_\2}(\Omega)),\\
  & \phi\partial_{t}(\rho_\alpha(p_\alpha)s_\alpha) \in
  L^2(0,T;(H^1_{\Gamma_\2}(\Omega))^\prime) \in
  L^2(0,T;(H^1_{\Gamma_\2}(\Omega))^\prime),
\end{align}
such that for all $\varphi,\, \psi \in
\mathbb{C}^1([0,T];H^1_{\Gamma_\2}(\Omega))\, \text{with} \,
\varphi(T,\cdot)=\psi(T,\cdot)=0$,
\begin{align}
  &-\int_{Q_T} \phi \rho_\2(p_\2) s_\2 \partial_t\varphi\dd x \dd t
  -\int_\Omega \phi (x) \rho_\2(p_\2^0(x)) s_\2^0(x)\varphi(0,x)\,dx \notag\\
  &+\int_{Q_T}M_\2(s_\2)\rho_\2(p_\2) { \nabla} p_\2 \cdot\nabla
  \varphi \dd x \dd t
  -\int_{Q_T}M_\2(s_\2)\rho_\2^2(p_\2){\bf g} \cdot\nabla \varphi \dd x \dd t\label{eq:p}\\
  &+\int_{Q_T}\rho_\2(p_\2)s_\2 f_{P}^{~}\varphi \dd x\dd t\notag =
    \int_{Q_T}\rho_\2(p_\2) s^I_\2 f_{I}^{~}\varphi \dd x \dd t,
\end{align}
\begin{align}
  &-\int_{Q_T} \phi \rho_\1(p_\1) s_\1 \partial_t\psi\dd x \dd t
  -\int_\Omega \phi (x) \rho_\1(p_\1^0(x)) s_\1^0(x)\psi(0,x)\,dx \notag\\
  &+\int_{Q_T}M_\1(s_\1)\rho_\1(p_\1) { \nabla} p_\1 \cdot\nabla \psi
  \dd x \dd t
  -\int_{Q_T}M_\1(s_\1)\rho_\1^2(p_\1){\bf g} \cdot\nabla \psi \dd x \dd t\label{eq:p}\\
  &+\int_{Q_T}\rho_\1(p_\1)s_\1 f_{P}^{~}\psi \dd x \dd t\notag
  =\int_{Q_T}\rho_\1(p_\1) s^I_\1 f_{I}^{~}\psi \dd x \dd t.
\end{align}
\end{definition}

\section{The finite volume scheme}
\label{sec:FV}
\subsection{Finite volume definitions and notations}
\label{subsec:FV}
\begin{figure}[ht]
 \begin{center}
   \includegraphics[width=0.4\linewidth]{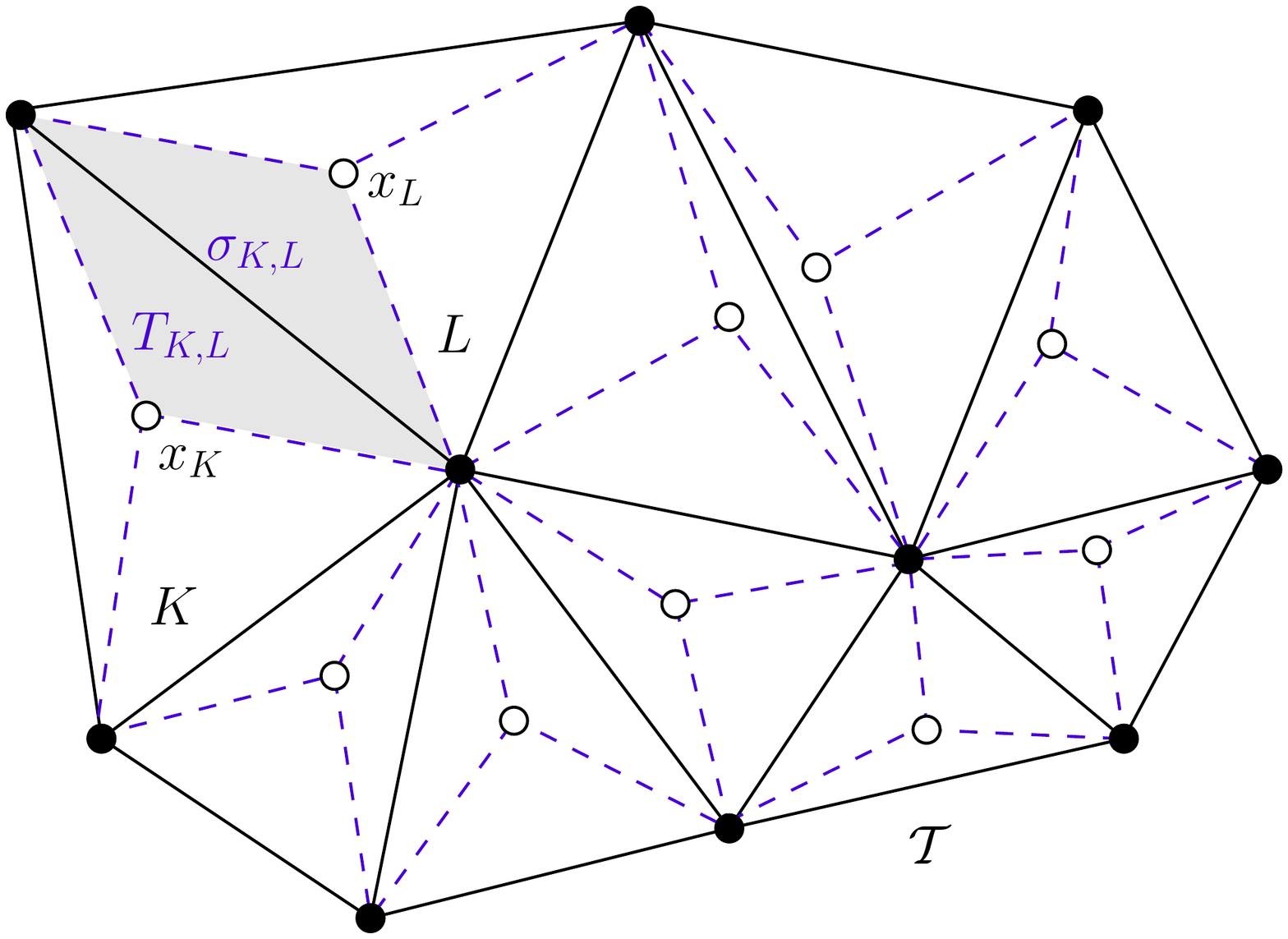}
 \end {center}
 \caption{\footnotesize
 Control volumes, centers and diamonds}
 \label{fig:diamond}
 \end{figure}

Following \cite{Eymard:book}, let us define a finite volume
discretization of $\Omega\times (0,T)$.
\begin{definition}\label{def:admissible_mesh} $\left( \text{Admissible mesh of } \Omega \right)$. An admissible mesh
  $\mathcal{T}$ of $\Omega$ is given by a set of open bounded
  polygonal convex subsets of $\Omega$ called control volumes and a
  family of points (the ``centers'' of control volumes) satisfying the
  following properties:
  \begin{enumerate}
  \item The closure of the union of all control volumes is
    $\overline{\Omega}$. We denote by $|K|$ the measure of $K$, and
    define 
    $$\sizt=\size(\mathcal{T}) = max\{diam(K),K\in \mathcal{T}
    \}.
    $$
  \item For any $(K,L)\in \T^2$ with $K \ne L$, then $K\cap L =
    \emptyset$. One denotes by $\E\subset \T^2$ the set of $(K,L)$
    such that the $d-1$-Lebesgue measure of $\overline{K}\cap
    \overline{L}$ is positive. For $(K,L)\in \E$, one denotes
    $\sigma_{\KL}=\overline{K}\cap\overline{L}$ and $|\sigma_{\KL}|$
    the $d-1$-Lebesgue measure of $\sigma_{\KL}$. And one denotes
    $\eta_{\KL}$ the unit normal vector to $\sig$ outward to $K$
  \item For any $K\in \T$, one defines $N(K) = \{L\in \T,
    (K,L)\in\E\}$ and one assumes that $\partial K = \overline{K}
    \backslash K = (\overline{K}\cap \partial\Omega)\cup (\cup_{L\in
      N(K)}\sigma_{\KL})$.
  \item The family of points $(x_K)_{K\in \T}$ is such that $x_K\in K$
    $(\text{for all }K\in \T)$ and, if $L\in N(K)$, it is assumed that
    the straight line $(x_K,x_L)$ is orthogonal to $\sigma_{\KL}$. We
    set $\dis =d(x_K,x_L)$ the distance between the points $x_K$ and $x_L$, and $\tokl=\frac {|\sigma_{\KL}|}
    {\dis}$, that is sometimes called the "transmissivity'' through
    $\sig$ (see Figure \ref{fig:diamond}).
    \item Let $\xi>0$. We assume the following regularity of the mesh : 
    \begin{equation}\label{eq:regularity}
    \forall K\in \T, \sum_{L\in N(K)} |\sig| \dis \le \xi |K|
  \end{equation}

\end{enumerate}
\end{definition}

 We denote by $H_\sizt(\Omega) \subset L^2(\Omega)$ the space of
 functions which are piecewise constant on each control volume $K \in
 \mathcal{T}$. For all $u_\sizt \in H_\sizt(\Omega)$ and for all $K
 \in \mathcal{T}$, we denote by $u_K$ the constant value of $u_\sizt$
 in $K$.  For $(u_\sizt,v_\sizt)\in (H_\sizt(\Omega))^2$, we define
 the following inner product:
$$
\left \langle u_\sizt,v_\sizt \right \rangle_{H_\sizt}= \dm\sum_{K \in
  \mathcal{T} }\sum_{L \in N(K) }
\frac{\abs{\sigma_{\KL}}}{\dis }(u_{L}-u_{K})(v_{L}-v_K),
$$
and the norm in $H_\sizt(\Omega)$ by
$$
\norm{u_\sizt}_{H_\sizt(\Omega)}=(\left \langle u_\sizt,u_\sizt \right
\rangle_{H_\sizt})^{1/2}.
$$
Finally, we define $L_\sizt(\Omega) \subset L^2(\Omega)$ the space of
functions which are piecewise constant on each control volume $K \in
\mathcal{T}$ with the associated norm
$$
\left (u_\sizt,v_\sizt \right )_{L_\sizt(\Omega)}=\sum_{K\in\mathcal{T}}\abs{K}u_{K} v_K,
\qquad
\norm{u_\sizt}^2_{L_\sizt(\Omega)} =\sum_{K \in \mathcal{T}}\abs{K}\abs{u_{K}}^2,
$$
for $(u_\sizt,v_\sizt)\in (L_\sizt(\Omega))^2$.
Further, a diamond $T_{\KL}$ is constructed upon the interface
$\sigma_{\KL}$, having $x_K$, $x_L$ for vertices (see Figure
\ref{fig:diamond}) and the $\dm$-dimensional mesure $\abs{T_{\KL}}$ of
$T_{\KL}$ equals to $\frac 1 \dm \abs{\sigma_{\KL}}\dis $.

The discrete gradient $\nabla_\sizt u_\sizt$ of a constant per control
volume function $u_\sizt$ is defined as the constant per diamond
$T_{\KL}$ $\R^\dm$-valued function with values
$$
\Grad_\sizt u_\sizt(x)=\begin{cases}
  & \dm \frac{u_{L}-u_{K}}{\dis }\eta_{\KL}\text{ if $x \in T_{\KL}$},\\
  & \dm \frac{u_{\sigma}-u_{K}}{d_{K,\sigma}}\eta_{K|\sigma}\text{ if
    $x \in T^{\text{ext}}_{K|\sigma}$}.
\end{cases}
$$
And the semi-norm $\|u_\sizt\|_{H_\sizt}$ coincides with the
$L^2(\Omega)$ norm of $\nabla_\sizt u_{\sizt}$, in fact
\begin{multline*}
  \|\Grad_\sizt u_\sizt\|_{L^2(\Omega)}^2 = \sum_{K \in \mathcal{T}}
  \sum_{L \in N(K) }\int_{T_{\KL}}|\nabla_\sizt u_\sizt|^2\,dx = \dm^2
  \sum_{K \in \mathcal{T} }\sum_{L \in N(K) }\abs{T_{\KL}}
  \frac{|u_{L}-u_{K}|^2}{|\dis |^2} \\
  = \dm \sum_{K \in \mathcal{T} }\sum_{L \in N(K)
  }\frac{\abs{\sigma_{\KL}} }{\dis }|u_{L}-u_{K}|^2 :=
  \norm{u_\sizt}_{H_\sizt(\Omega)}^2.
 \end{multline*}
 We assimilate a discrete field $(\vec{F}_{\KL})$ on $\Omega$ to the
 piecewise constant vector-function 
 $$
 \vec{F}_\sizt = \sum_{\sigma_{\KL}\in\mathcal{E}} \vec{F}_{\KL}
 \charu_{T_{\KL}}.
 $$ 
 The discrete divergence of the field $\vec{F}_\sizt$ is defined as
 the discrete function $w_\sizt=div_\sizt\vec{F}_\sizt $ with the
 entires
 \begin{equation}\label{equality:div}   
   div_K \vec{F}_\sizt :=\frac{1}{\abs{K}} \sum_{L\in N(K)}
   |\sigma_{\KL}|\vec{F}_{\KL}\cdot \eta_{\KL}.
\end{equation}

\bigskip

The problem under consideration is time-dependent, hence we also need
to discretize the time interval $(0,T)$.
\begin{definition}\label{def:disc_time} $\left( \text{Time discretization}\right)$. A time discretization of
  $(0,T)$ is given by an integer value $N$ and by a strictly
  increasing sequence of real values $(t^n)_{n\in [0,N+1]}$ with
  $t^0=0$ and $t^{N+1} = T$. Without restriction, we consider a uniform step time  $\delta t = t^{n+1} - t^{n}$, for $n\in [0,N]$.
\end{definition}

We may then define a discretization of the whole domain $\Omega \times (0,T)$ in the following way:
\begin{definition}\label{def:disc_QT}$\left( \text{Discretization of }\Omega \times (0,T)\right)$. A finite volume discretization $\D$ of $\Omega \times (0,T)$ is defined by
  $$
  \D = \Big(\T,\E,(x_K)_{K\in\T},N,(t^n)_{n\in[0,N]}\Big),
  $$
  where $\T,\E,(x_K)_{K\in\T}$ is an admissible mesh of $\Omega$ in
  the sense of Definition \ref{def:admissible_mesh} and
  $N,(t^n)_{n\in[0,N]}$ is a time discretization of $(0,T)$ in the
    sense of Definition \ref{def:disc_time}. One then sets 
    $$
    \size(\D) = max(\size(\T),\delta t). 
    $$
   \end{definition}

\begin{definition}\label{def:disc_function}$\left(\text{Discrete functions and
      notations}\right)$. Let $\D$ be a discretization of
  $\Omega\times (0,T)$ in the sense of Definition \ref{def:disc_QT}.
  We denote any function from $\T\times [0,N+1]$ to $\R$ by using the
  subscript $\D$, $(s_{\alpha,\D} \text{ and } p_{\alpha,\D} \text{
    for instance})$ and we denote its value at the point $(x_K,t^n)$ using
  the subscript $K$ and the superscript $n$ $(s_{\alpha,K}^n \text{
    for instance, we then denote }
  s_{\alpha,\D}=(s_{\alpha,K}^n)_{K\in\T,n\in[0,N+1]} )$. To any
  discrete function $u_{\D}$ corresponds an approximate function
  defined almost everywhere on $\Omega\times(0,T)$ by:
  $$
  u_{\D}(t,x) = u_{K}^{n+1}, \text{ for a.e. }(t,x)\in
  (t^n,t^{n+1})\times K, \forall K\in \T, \forall n\in [0,N].
  $$  
  For any continuous function $f:\R\to\R$, $f(u_{\D})$ denotes the
  discrete function $(K,n)\to f(u_{K}^{n+1})$. if $L\in N(K)$, and
  $u_{\D}$ is a discrete function, we denote by $\delta_{\KL}^{n+1}(u)
  = u_{L}^{n+1}-u_{K}^{n+1}$. For example, $\de(f(u)) = f(u_{L}^{n+1})
  - f(u_{K}^{n+1})$.
 \end{definition}

\bigskip

Let us recall the following two lemmas :
\begin{lemma}\label{lem:poincare}$\left(\text{Discrete Poincar\'{e}
      inequality } \right)$ \cite{Eymard:book}. Let $\Omega$ be an
  open bounded polygonal subset of $\R^\dm$, $ \dm=2~ \text{or}~ 3 $.
  Let $\T$ be a finite volume discretization of $\Omega$ in the sense
  of Definition \ref{def:admissible_mesh}, and let $u$ be a function
  which is constant on each cell $K \in \mathcal{T}$, that is,
  $u(x)=u_K$ if $x\in K,$ then
$$\norm{u}_{L^2(\Omega)}\leq diam(\Omega) \norm{u}_{H_h(\Omega)},$$
where $\norm{\cdot}_{H_h(\Omega)}$ is the discrete $H_0^1$ norm.
\end{lemma}
\begin{remark}$\left(\text{Dirichlet condition on part of the
      boundary} \right)$. The \it{lemma} \ref{lem:poincare} gives a
  discrete Poincar\'{e} inequality for Dirichlet boundary conditions
  on the boundary $\partial \Omega$. In the case of Dirichlet
  condition on part of the boundary only, it is still possible to
  prove a discrete Poincar\'{e} inequality provided that the polygonal
  bounded open set $\Omega$ is connected.
\end{remark}
\begin{lemma}\label{lem:integration_parti}$\left(\text{Discrete integration 
      by parts formula}\right)$.  
  Let $F_{K/L},~ K\in
  \mathcal{T}$ and $L\in N(K)$ be a value in $\R$ depends on $K$ and
  $L$ such that $F_{K/L}=-F_{L/K}$ and let $\varphi$ be a function
  which is constant on each cell $K \in \mathcal{T}$, that is,
  $\varphi(x)=\varphi_K$ if $x\in K,$ then
\begin{equation}
\label{integration_part}
\sum_{K\in \mathcal{T}}\sum_{L\in N(K)}F_{K/L} \varphi_K
=-\frac{1}{2}\sum_{K\in \mathcal{T}}\sum_{L\in N(K)}
F_{K/L}(\varphi_L-\varphi_K)
\end{equation}
Consequently, if $F_{K/L}= a_{K/L}(b_L-b_K)$, with $a_{K/L}=a_{L/K}$,
then
\begin{equation}
  \label{part}
  \sum_{K\in \mathcal{T}}\sum_{L\in
    N(K)}a_{K/L}(b_L-b_K)\varphi_K=-\frac{1}{2}\sum_{K\in
    \mathcal{T}}\sum_{L\in N(K)}a_{K/L}(b_L-b_K
  )(\varphi_L-\varphi_K)
\end{equation}
\end{lemma}
%

 \subsection{The coupled finite volume scheme}
 The finite volume scheme is obtained by writing the balance equations
 of the fluxes on each control volume. Let $\mathcal{\D}$ be a
 discretization of $\Omega\times (0,T)$ in the sense of Definition
 \ref{def:disc_QT}. Let us integrate equations
 \eqref{eq:conservation_masse} over each control volume $K$. By using
 the Green formula, if $\Phi$ is a vector field, the integral of $\di
 (\Phi)$ on a control volume $K$ is equal to the sum of the normal
 fluxes of $\Phi$ on the edges \eqref{equality:div}. Here we apply
 this formula to approximate $M_\alpha(s_\alpha)\Grad
 p_\alpha\cdot\eta_{\KL},~~ (\alpha=\2,\1)$ by means of the values
 $s_{\alpha,K},s_{\alpha,L}$ and $p_{\alpha,K},p_{\alpha,L}$ that are
 available in the neighborhood of the interface $\sigma_{\KL}$. To do
 this, let us use some function $G_\alpha$ of $(a,b,c)\in \R^3$ . The
 numerical convection flux functions $G_\alpha \in C(\R^3,\R)$, are
 required to satisfy the properties:
\begin{equation}\label{Hypfluxes}
\begin{cases}
  \text{(a) $G_\alpha(\cdot,b,c)$
    is non-decreasing for all $b,c\in\R$,}\\
  \hspace*{15pt} \text{and $G_\alpha(a,\cdot,c)$ is non-increasing
    for all $ a,c\in \R$};\\
  \text{(b) $G_\alpha(a,a,c)=-M_\alpha(a)\,c$ for all $a,c \in \R$}; \\
  \text{(c) $G_\alpha(a,b,c)=-G_\alpha(b,a,-c)$
    and there exists $C>0$ such that }\\
  \,\,\quad|G_\alpha(a,b,c)|\leq C \, \bigl(|a|+|b|\bigr)|c| \text{
    for all $a,b,c \in \R$}.
\end{cases}
\end{equation}

Note that the assumptions (a), (b) and (c) are standard and they
respectively ensure the maximum principle on saturation, the
consistency of the numerical flux and the conservation of the
numerical flux on each interface. Practical examples of numerical
convective flux functions can be found in \cite{Eymard:book}.

In our context, we consider an upwind scheme, the numerical flux $G_\alpha$ satisfying
  \eqref{Hypfluxes} defined by
\begin{align}\label{Galphaupwind}
  G_\alpha(a,b,c)=-M_\alpha(b)\,{c}^+ +M_\alpha(a) \,{c}^-
\end{align}
where $c^+ = \max(c,0)$ and $c^- = \max(-c,0)$. Note that the function
$s_\alpha \mapsto M_\alpha(s_\alpha)$ is non-decreasing, which lead to
the monotony property of the function $G_\alpha$.

\bigskip

The resulting equation is discretized with a  implicit Euler scheme in time;
 the normal gradients are discretized with a
centered finite difference scheme.\\
Denote by
$p_{\alpha,\D}=(p_{\alpha,K}^{n+1})_{K \in \mathcal{T},n\in[0,N]}$ and
$s_{\alpha,\D}=(s_{\alpha,K}^{n+1})_{K \in \mathcal{T},n\in[0,N]}$ the
discrete unknowns corresponding to $p_\alpha$ and $s_\alpha$.
 The finite volume scheme is the following set of equations :
\begin{equation}\label{eq:p0} 
  p_{\alpha,K}^0=\frac{1}{\abs{K}} \int_{K} p^0_\alpha(x) \dd x,\;  
  s_{\alpha,K}^0=\frac{1}{\abs{K}} \int_{K} s^0_\alpha(x) \dd x, \text{ for all } K\in\T,
\end{equation} 

\begin{multline}\label{sys_disc:pl}
  \abs{K}\phi_K\frac{\rho_\2(p^{n+1}_{\2,K})s^{n+1}_{\2,K}-\rho_\2(p^{n}_{\2,K})s^{n}_{\2,K}}{\h}
  + \sum_{L \in N(K) }\tokl
  \rho^{n+1}_{\2,\KL}G_\2(s^{n+1}_{\2,K},s^{n+1}_{\2,L};\de(p_\2)) \\
  +F^{\;n+1}_{\2,K} +\abs{K} \rho_\2(p_{\2,K}^{n+1}) s_{\2,K}^{n+1}
  f_{P,K}^{n+1} = \abs{K} \rho_\2(p_{\2,K}^{n+1}) \sli f_{I,K}^{n+1},
\end{multline}
\begin{multline}\label{sys_disc:pg}
  \abs{K}\phi_K\frac{\rho_\1(p^{n+1}_{\1,K})s^{n+1}_{\1,K}-\rho_\1(p^{n}_{\1,K})s^{n}_{\1,K}}{\h}
  + \sum_{L \in N(K) }\tokl
  \rho^{n+1}_{\1,\KL}G_\1(s^{n+1}_{\1,K},s^{n+1}_{\1,L};\de(p_\1)) \\
  +F^{\;n+1}_{\1,K} +\abs{K} \rho_\1(p_{\1,K}^{n+1}) s_{\1,K}^{n+1}
  f_{P,K}^{n+1} = \abs{K} \rho_\1(p_{\1,K}^{n+1}) \sgi f_{I,K}^{n+1},
\end{multline}
\begin{equation}\label{def_disc:pc}
  p_c(s^{n+1}_{\1,K}) = p_{\2,K}^{n+1}-p_{\1,K}^{n+1},
\end{equation}
where $\disp F^{n+1}_{\alpha,K}$ $(\alpha=\2,~\1)$ the approximation
of $\disp \int_{\partial K}\rho_\alpha^2(p_\alpha^{n+1})
M_\alpha(s_\alpha^{n+1}) {\bf g}\cdot \eta_{\KL}\, \dd \Gamma(x)$ by
an upwind scheme:
\begin{multline} \label{fgravity}
  F^{n+1}_{\alpha,K}=\sum_{L \in N(K) }F^{n+1}_{\alpha,\KL}=\sum_{L
    \in N(K) }|
  \sig|(\rho^{n+1}_{\alpha,\KL})^2 \Big(M_\alpha(s^{n+1}_{\alpha,K}){\bf g}_{\KL}^+
  -  M_\alpha(s^{n+1}_{\alpha,L}){\bf g}_{\KL}^-\Big),
\end{multline}
with  ${\bf g}_{\KL}^+ :=    ({\bf g}\cdot \eta_{\KL})^+ $ and  ${\bf g}_{\KL}^- :=  ({\bf g}\cdot \eta_{\KL})^- $.
Notice that  the source terms are, for  $n\in \{0,\ldots,N-1\}$
$$
f^{n+1}_{P,K}:=
          \frac{1}{\h\abs{K}}
          \int_{t^n}^{t^{n+1}}
          \int_K f_P(t,x)\,dx dt,\quad f^{n+1}_{I,K}:=
          \frac{1}{\h\abs{K}}
          \int_{t^n}^{t^{n+1}}
          \int_K f_I(t,x)\,dx dt 
         $$
\bigskip

The mean value of the density of each phase on interfaces is not
classical since it is given as
\begin{equation}\label{meanrho}
\begin{aligned}
  \frac{1}{\rho^{n+1}_{\alpha,\KL}}=\begin{cases}
    \frac{1}{p_{\alpha,L}^{n+1}-p_{\alpha,K}^{n+1}}
    \int_{p_{\alpha,K}^{n+1}}^{p_{\alpha,L}^{n+1}}
    \frac{1}{\rho_\alpha(\zeta)}\,d\zeta
    & \text{ if } p_{\alpha,K}^{n+1} \ne p_{\alpha,L}^{n+1}, \\
    \frac{1}{\rho^{n+1}_{\alpha,K}} & \text{ otherwise}.
          \end{cases}
\end{aligned}
 \end{equation}
 This choice is crucial to obtain estimates on discrete pressures.

Note that the numerical fluxes to approach the gravity terms
$F_\alpha$ are nondecreasing with respect to $s_{\alpha,K}$ and
nonincreasing with respect to $s_{\alpha,L}$.

The upwind fluxes \eqref{Galphaupwind} can be rewritten in the equivalent form 
  \begin{equation}
  \label{Gupwind}
  G_\alpha(s^{n+1}_{\alpha,K},s^{n+1}_{\alpha,L};
  \de(p_{\alpha}))=-M_{\alpha}(s_{\alpha,K|L}^{n+1})\;\de(p_{\alpha}),
\end{equation}
where $M_{\alpha}(s_{\alpha,K|L}^{n+1})$ denote the upwind
discretization of $M_\alpha(s_\alpha)$ on the interface $\sigma_{\KL}$
and
\begin{align}\label{notation:saturation_interface}
s_{\alpha,K|L}^{n+1}=\begin{cases}
  & s_{\alpha,K}^{n+1} \text{ if } (K,L)\in \mathcal{E}_\alpha^{n+1},\\
  & s_{\alpha,L}^{n+1} \text{ otherwise, }
\end{cases}
\end{align}
with the set $\mathcal{E}_\alpha^{n+1}$ is subset of $\mathcal{E}$
such that
\begin{align}\label{notation:set}
  \mathcal{E}_\alpha^{n+1}=\{ (K,L)\in \mathcal{E},
  \de(p_\alpha)=p_{\alpha,L}^{n+1}-p_{\alpha,K}^{n+1}\le 0 \}.
\end{align}
We extend the mobility functions $s_\alpha\mapsto M_\alpha(s_\alpha)$
outside $[0,1]$ by continuous constant functions. We show below (see
Prop. \ref{prop:existance}) that there exists at least one solution to
this scheme. From this discrete solution, we build an approximation
solution $p_{\alpha,\D}$ defined almost everywhere on $Q_T$ by (see
Definition \ref{def:disc_function}):
\begin{equation}\label{prob:general}
  p_{\alpha,\D}(t,x)=p_{\alpha,K}^{n+1},\;\forall x \in K,\forall t \in (t^n,t^{n+1}).
\end{equation}
The main result of this paper is the following theorem.
\begin{theorem} \label{theo:principal} Assume hypothesis
  (H\ref{hyp:H1})-(H\ref{hyp:H6}) hold. Let $\{\D_m\}_{m\in \N}$ be a
  sequence of discretization of $Q_T$ in the sense of definition
  \ref{def:disc_QT} such that $\lim_{m\to +\infty}\size(\D_m)=0$. Let
  $(p^{0}_\alpha,s_\alpha^0) \in L^2(\Omega, \R)\times
  L^\infty(\Omega, \R)$. Then there exists an approximate solutions
  $(p_{\alpha,\D_m})_{m\in \N}$ corresponding to the system
  (\ref{sys_disc:pl})-(\ref{sys_disc:pg}), which converges (up to a
  subsequence) to a weak solution $p_\alpha$ of
  \eqref{eq:conservation_masse} in the sense of the Definition
  \ref{def}.
\end{theorem}

\section{Preliminary fundamental lemmas}\label{sec:fundamental_lemma}
The mobility of each phase vanishes in the region where
the phase is missing. Therefore, if we control the quantities
$M_\alpha \nabla p_\alpha$ in the $L^2$-norm, this does not permit the
control of the gradient of pressure of each phase. In the continuous
case, we have the following relationship between the global pressure,
capillary pressure and the pressure of each phase
\begin{equation}\label{eq:relation}
  \begin{aligned}
    M |\nabla p |^{2} + \frac{M_{l}M_{g}}{M}|\nabla p_c |^{2} = M_{l}
    |\nabla p_{l}|^2 + M_{g} |\nabla p_{g}|^2.
  \end{aligned}
\end{equation}
This relationship, means that, the control of the velocities ensures
the control of the global pressure and the capillary terms $\B$ in the
whole domain regardless of the presence or the disappearance of the
phases. This estimates (of the global pressure and the capillary terms
$\B$) has a major role in the analysis, to treat the degeneracy of the
dissipative terms $\di(\rho_\alpha M_\alpha \nabla p_\alpha)$.

In the discrete case, these relationship, are not obtained in a
straightforward way. This equality is replaced by three discrete
inequalities which we state in the following three lemmas.  

We derive in the next lemma the preliminary step to proof the
estimates of the global pressure and the capillary terms given in
Proposition \ref{prop:estimation_pression} and Corollary
\ref{cor:est}. These lemmas are first used to prove a compactness lemma
and then used for the convergence result.

\begin{lemma}\label{lemma:preliminary_step}(Total mobility and global pressure \cite{Eymard00}).
  Under the assumptions $({H}\ref{hyp:H1})-({H}\ref{hyp:H6})$ and the
  notations \eqref{def:pression_globale}. Let $\D$ be a finite volume
  discretization of $\Omega\times (0,T)$ in the sense of Definition
  \ref{def:disc_QT}.
  Then for all $(K,L)\in \mathcal{E}$ and for all $n\in [0,N]$ the
  following inequalities hold:
  \begin{align}\label{est:mobilite}
    M_{\2,K|L}^{n+1} + M_{\1,K|L}^{n+1} \ge m_0,
  \end{align}
\text{ and }
 \begin{align}\label{ineq:pression_globale}
   m_0\Big(\de (p)\Big)^2
   \le M_{\2,K|L}^{n+1}\Big(\de
   (p_\2)\Big)^2 + M_{\1,K|L}^{n+1}\Big(\de (p_\1)\Big)^2.
 \end{align}
\end{lemma}

The proof of this lemma is made  by R. Eymard and al. in \cite{Eymard00}. The proof of this result can be applied for compressible flow since the proof use only the definition of the global pressure.

\begin{lemma}\label{prellem2}$\left(\text{Capillary term } \B \right)$.
  Under the assumptions $({H}\ref{hyp:H1})-({H}\ref{hyp:H6})$ and the
  notations \eqref{def:pression_globale}. Let $\D$ be a finite volume
  discretization of $\Omega\times (0,T)$ in the sense of Definition
  \ref{def:disc_QT}.
  Then there exists a constant $C>0$ such that for all $(K,L)\in  \mathcal{E}$ and $ n\in [0,N]$ : 
  \begin{equation}\label{ineq:discrete_beta}
    (\de(\B(s_\3)))^2  \le M_{\2,K|L}^{n+1}\Big(\de (p_\2)\Big)^2 +
    M_{\1,K|L}^{n+1}\Big(\de (p_\1)\Big)^2. 
  \end{equation}
\end{lemma}
In the incompressible case (see \cite{Eymard00}) this kind of estimate
is obtained by using the mass conservation equation and under
hypotheses ont the relative permeability of the $\alpha$ phase,
whereas, the compressibility add more difficulties, our approach use
only the definition of the function $\B$ and consequently this lemma
can be used for compressible and incompressible degenerate flows.
\begin{proof}  We take the same
  decomposition of the interface as that proposed by R. Eymard and al.
  in \cite{Eymard00}, namely the different possible cases
  $(K,L)\in\E_\2^{n+1}\cap\E_\1^{n+1}$,
  $(K,L)\notin\E_\2^{n+1}\cup\E_\1^{n+1}$, $(K,L)\in\E_\2^{n+1}$ and
  $(K,L)\notin\E_\1^{n+1}$, and the last case $(K,L)\notin\E_\2^{n+1}$
  and $(K,L)\in\E_\1^{n+1}$; where the sets $\E_\2^{n+1}$ and $\E_\1^{n+1}$ are defined in \eqref{notation:set}. We establish \label{ineq:discrete_beta} for the four cases. 
  
    $\bullet${\bf First case.} If $(K,L)\notin \E_\2$ and $(K,L)\in
  \E_\1$.  We may notice that if the upwind choice is different for
  the two equations, we have
$$
M_{\alpha,K|L}^{n+1} = \max_{[s_{\3,K},s_{\3,L}]} M_\alpha.
$$
By definition of $\B$ in \eqref{def:beta}, there exists some $a\in
[s_{\3,K},s_{\3,L}]$ such that 
$$\de(\B(s_\3)) = -
\frac{M_\2(a)M_\1(a)}{M_\2(a)+M_\1(a)}\de(p_c(s_\3)),$$
 we then get 
\begin{align*}
  (\de(\B(s_\3)))^2 & \le M_{\2,K|L}^{n+1}M_{\1,K|L}^{n+1}
  (\de(p_c(s_\3)))^2 \\ & \le C_1 M_{\2,K|L}^{n+1}M_{\1,K|L}^{n+1}
  \Big((\de(p_\1))^2+(\de(p_\2))^2 \Big) \\ & \le C_2
  \Big(M_{\2,K|L}^{n+1}(\de(p_\2))^2+M_{\1,K|L}^{n+1}(\de(p_\1))^2
  \Big).
\end{align*}
\\
$\bullet${\bf Second case:} The case $(K,L)\in \E_\2$ and
$(K,L)\notin \E_\1$ is similar.\\
$\bullet${\bf Third case:} The case $(K,L) \in \E_\2$ and $(K,L)\in
\E_\1$.  We have
\begin{equation}\label{relation:ship}
\begin{aligned}
  & M_{\2,K|L}^{n+1}(\de(p_\2))^2+M_{\1,K|L}^{n+1}(\de(p_\1))^2 \\ &=
  M_\2(s_{\2,K}^{n+1})(\de(p_\2))^2+M_\1(s_{\1,K}^{n+1})(\de(p_\1))^2
  \\ &= \Big(M_\2(s_{\2,K}^{n+1}) + M_\1(s_{\1,K}^{n+1})\Big)
  (\de(p))^2 \\ & + M_\1(s_{\1,K}^{n+1})(\de(\tilde{p}(s_\3)))^2 +
  M_\2(s_{\2,K}^{n+1})(\de(\bar{p}(s_\3)))^2 \\ & - 2
  M_\1(s_{\1,K}^{n+1})\de(p) \de(\tilde{p}(s_\3)) -2
  M_\2(s_{\2,K}^{n+1})\de(p)\de(\bar{p}(s_\3)).
\end{aligned}
\end{equation}
We will distinguish the case $s_{\3,K}^{n+1}\le s_{\3,L}^{n+1}$ and the case $s_{\3,K}^{n+1}\ge s_{\3,L}^{n+1}$.
\begin{enumerate}
\item If we assume that $s_{\3,K}^{n+1}\le s_{\3,L}^{n+1}$, we deduce that 
  \begin{enumerate}
  \item $\de(\bar{p}(s_\3)) \le 0 \text{ since } \bar {p}(s_\3) \text{
      is nonincreasing}$,
  \item $\de(\tilde{p}(s_\3)) \ge 0 \text{ since } \tilde {p}(s_\3)
    \text{ is nondecreasing}$,
  \item $\de(p)=\de(p_\2)+\de(\bar{p}(s_\3))\le 0$.
  \end{enumerate}
  One then gets from \eqref{relation:ship} that:
\begin{align*}
  & M_\2(s_{\2,K}^{n+1})(\de(p_\2))^2+M_\1(s_{\1,K}^{n+1})(\de(p_\1))^2 \\
  &\ge \Big(M_\2(s_{\2,K}^{n+1}) + M_\1(s_{\1,K}^{n+1})\Big) (\de(p))^2 \\
  & + M_\1(s_{\1,K}^{n+1})(\de(\tilde{p}(s_\3))^2) +
  M_\2(s_{\2,K}^{n+1})(\de(\bar{p}(s_\3)))^2 \\ & - 2
  M_\2(s_{\2,K}^{n+1} )\de(p) \de(\bar{p}(s_\3)).
\end{align*}
The previous inequality gives:
\begin{align*}
  & \Big(M_\2(s_{\2,K}^{n+1}) + M_\1(s_{\1,K}^{n+1})\Big) (\de(p))^2 \\
  & + M_\1(s_{\1,K}^{n+1})(\de(\tilde{p}(s_\3)))^2 +
  M_\2(s_{\2,K}^{n+1})(\de(\bar{p}(s_\3)))^2 \\ & \le
  M_\1(s_{\1,K}^{n+1})(\de(p_\1))^2+M_\2(s_{\2,K}^{n+1})(\de(p_\2))^2+2
  M_\2(s_{\2,K}^{n+1} )\de(p) \de(\bar{p}(s_\3))\\ & \le
  M_\2(s_{\2,K}^{n+1})(\de(p_\2))^2+M_\1(s_{\1,K}^{n+1})(\de(p_\1))^2
  \\ & + M_\2(s_{\2,K}^{n+1}) (\de(p))^2 + M_\2(s_{\2,K}^{n+1})(
  \de(\bar{p}(s_\3)))^2,
\end{align*}
which implies the inequality:
\begin{equation}
\begin{aligned}
  & M_\1(s_{\1,K}^{n+1})(\de(p))^2 +
  M_\1(s_{\1,K}^{n+1})(\de(\tilde{p}(s_\3)))^2 \\ & \le
  M_\2(s_{\2,K}^{n+1})(\de(p_\2))^2+M_\1(s_{\1,K}^{n+1})(\de(p_\1))^2.
\end{aligned}
\end{equation}
Or, by definition of $\B$ \eqref{def:beta}, there exists some $a\in
[s_{\3,K},s_{\3,L}]$ such that $\de(\B(s_\3)) =M_\1(a)
\de(\tilde{p}(s_\3))$, we get then
 \begin{align*}
   (\de(\B(s_\3)))^2 & \le M_\1(s_{\1,K}) (\de(\tilde{p}(s_\3)))^2\\ &
   \le
   M_\2(s_{\2,K}^{n+1})(\de(p_\2))^2+M_\1(s_{\1,K}^{n+1})(\de(p_\1))^2.
  \end{align*}
  which is \eqref{ineq:discrete_beta} in that case. \\
  %
\item If we assume that $s_{\3,L}^{n+1}\le s_{\3,K}^{n+1}$, we get
  that
  \begin{enumerate}
  \item $\de(\bar{p}(s_\3)) \ge 0 \text{ since } \bar {p}(s_\3) \text{
      is nonincreasing}$,
  \item $\de(\tilde{p}(s_\3)) \le 0 \text{ since } \tilde {p}(s_\3)
    \text{ is nondecreasing}$,
  \item $\de(p)=\de(p_\1)+\de(\tilde{p}(s_\3))\le 0$.
  \end{enumerate}
  One then gets from \eqref{relation:ship} that:
\begin{align*}
  & M_\2(s_{\2,K}^{n+1})(\de(p_\2))^2+M_\1(s_{\1,K}^{n+1})(\de(p_\1))^2 \\
  &\ge \Big(M_\1(s_{\1,K}^{n+1}) + M_\2(s_{\2,K}^{n+1})\Big) (\de(p))^2 \\
  & + M_\1(s_{\1,K}^{n+1})(\de(\tilde{p}(s_\3))^2) +
  M_\2(s_{\2,K}^{n+1})(\de(\bar{p}(s_\3)))^2 \\ & - 2
  M_\1(s_{\1,K}^{n+1} )\de(p) \de(\tilde{p}(s_\3)).
\end{align*}
The previous inequality gives:
\begin{align*}
  & \Big(M_\2(s_{\2,K}^{n+1}) + M_\1(s_{\1,K}^{n+1})\Big) (\de(p))^2 \\
  & + M_\1(s_{\1,K}^{n+1})(\de(\tilde{p}(s_\3))^2) +
  M_\2(s_{\2,K}^{n+1})(\de(\bar{p}(s_\3)))^2 \\ & \le
  M_\2(s_{\2,K}^{n+1})(\de(p_\2))^2+M_\1(s_{\1,K}^{n+1})(\de(p_\1))^2+2
  M_\1(s_{\1,K}^{n+1} )\de(p) \de(\tilde{p}(s_\3))\\ & \le
  M_\2(s_{\2,K}^{n+1})(\de(p_\2))^2+M_\1(s_{\1,K}^{n+1})(\de(p_\1))^2
  \\ & + M_\1(s_{\1,K}^{n+1}) (\de(p))^2 + M_\1(s_{\1,K}^{n+1})(
  \de(\tilde{p}(s_\3)))^2,
\end{align*}
which implies the inequality:
\begin{equation}
\begin{aligned}
  & M_\2(s_{\2,K}^{n+1})(\de(p))^2 +
  M_\2(s_{\2,K}^{n+1})(\de(\bar{p}(s_\3)))^2 \\ & \le
  M_\2(s_{\2,K}^{n+1})(\de(p_\2))^2+M_\1(s_{\1,K}^{n+1})(\de(p_\1))^2.
\end{aligned}
\end{equation}
Or, by definition of $\B$ \eqref{def:beta} there exists some $a\in
[s_{\3,K},s_{\3,L}]$ such that $\de(\B(s_\3)) = -M_\2(a)
\de(\bar{p}(s_\3))$, we get then
\begin{align*}
  (\de(\B(s_\3)))^2 & \le M_\2(s_{\2,K}) (\de(\bar{p}(s_\3)))^2\\ &
  \le
  M_\2(s_{\2,K}^{n+1})(\de(p_\2))^2+M_\1(s_{\1,K}^{n+1})(\de(p_\1))^2,
\end{align*}
\end{enumerate}
which is \eqref{ineq:discrete_beta} in that case.\\

$\bullet${\bf Fourth case:} The case $(K,L)\notin \E_\2$ and
$(K,L)\notin \E_\1$ is similar of the third case.
\end{proof}

\begin{lemma}$\left(\text{Dissipative terms}\right)$. Under the
  assumptions $({H}\ref{hyp:H1})-({H}\ref{hyp:H6})$ and the notations
  \eqref{def:pression_globale}. Let $\D$ be a finite volume
  discretization of $\Omega\times (0,T)$ in the sense of Definition
  \ref{def:disc_QT}.
  Then there exists a constant $C>0$ such that for all $(K,L)\in \mathcal{E}$ and $n\in
    [0,N]$
  \begin{equation}\label{ineq:discrete_pbar}
    M_{\2,K|L}^{n+1}(\de(\bar{p}(s_\3)))^2  \le
    M_{\2,K|L}^{n+1}\Big(\de (p_\2)\Big)^2 + M_{\1,K|L}^{n+1}\Big(\de
    (p_\1)\Big)^2,
      \end{equation}
and 
   \begin{equation}\label{ineq:discrete_ptilde}
     M_{\1,K|L}^{n+1}(\de(\tilde{p}(s_\3)))^2  \le
     M_{\2,K|L}^{n+1}\Big(\de (p_\2)\Big)^2 + M_{\1,K|L}^{n+1}\Big(\de
     (p_\1)\Big)^2.
   \end{equation}
\end{lemma}

\begin{proof} In order to prove \eqref{ineq:discrete_pbar} and
  \eqref{ineq:discrete_ptilde}, we consider the exclusive cases
  $(K,L)\in \E_\2^{n+1}\cap\E_\1^{n+1}$, $(K,L)\notin
  \E_\2^{n+1}\cup\E_\1^{n+1}$, $(K,L)\notin \E_\2^{n+1}$ and $(K,L)\in
  \E_\1^{n+1}$
  and the last case $(K,L)\in \E_\2^{n+1}$ and $(K,L)\notin \E_\1^{n+1}$.
  
  {\bf First case.} If $(K,L)\notin \E_\2$ and $(K,L)\in
  \E_\1$. We have
$$
M_{\alpha,K|L}^{n+1} = \max_{[s_{\3,K},s_{\3,L}]} M_\alpha,
$$
and by definition of $\bar{p}$ there exists some $a\in
[s_{\3,K},s_{\3,L}]$ such that $\de(\bar{p}(s_\3)) =
\frac{M_\1(a)}{M_\1(a)+M_\2(a)}\de(p_c(s_\3))$, we get then
\begin{align*}
  M_{\2,K|L}^{n+1}(\de(\bar{p}(s_\3)))^2 & \le
  M_{\2,K|L}^{n+1}M_{\1,K|L}^{n+1} (\de(p_c(s_\3)))^2 \\ & \le C_1
  M_{\2,K|L}^{n+1}M_{\1,K|L}^{n+1} \Big((\de(p_\1))^2+(\de(p_\2))^2
  \Big) \\ & \le C_2
  \Big(M_{\2,K|L}^{n+1}(\de(p_\2))^2+M_{\1,K|L}^{n+1}(\de(p_\1))^2
  \Big),
\end{align*}
which gives \eqref{ineq:discrete_pbar}. For the discrete estimate
\eqref{ineq:discrete_ptilde} and by definition of $\tilde{p}$ there
exists some $b\in [s_{\3,K},s_{\3,L}]$ such that $\de(\tilde{p}(s_\3))
= - \frac{M_\2(b)}{M_\1(b)+M_\2(b)}\de(p_c(s_\3))$, we get then
\begin{align*}
  M_{\1,K|L}^{n+1}(\de(\tilde{p}(s_\3)))^2 & \le
  M_{\1,K|L}^{n+1}M_{\2,K|L}^{n+1} (\de(p_c(s_\3)))^2 \\ & \le C_1
  M_{\1,K|L}^{n+1}M_{\2,K|L}^{n+1} \Big((\de(p_\1))^2+(\de(p_\2))^2
  \Big) \\ & \le C_2
  \Big(M_{\2,K|L}^{n+1}(\de(p_\2))^2+M_{\1,K|L}^{n+1}(\de(p_\1))^2
  \Big),
\end{align*}
which gives \eqref{ineq:discrete_ptilde}.

{\bf Second case.} The case $(K,L) \in \E_\2$ and
$(K,L)\notin E_\1$ is similar.

The {\bf third case} and the {\bf fourth case} can be treated as the cases in the lemma \ref{prellem2}.
\end{proof}

\section{A priori estimates and existence of the approximate
  solution}\label{sec:basic-apriori}

We derive new energy estimates on the discrete velocities
$M_\alpha(s_{\alpha,K|L}^{n+1}) \de(p_\alpha)$. Nevertheless, these
estimates are degenerate in the sense that they do not permit the
control of $\de(p_\alpha)$, especially when a phase is missing. So,
the global pressure has a major role in the analysis, we will show
that the control of the discrete velocities
$M_\alpha(s_{\alpha,K|L}^{n+1}) \de(p_\alpha)$ ensures the control of
the discrete gradient of the global pressure and the discrete gradient
of the capillary term ${\mathcal B}$ in the whole domain regardless of the
presence or the disappearance of the phases.

The following section gives us some necessary energy estimates to
prove the theorem \ref{theo:principal}.
\subsection{The maximum principle}
Let us show in the following Lemma that the phase by phase upstream
choice yields the $L^\infty$ stability of the scheme which is a basis
to the analysis that we are going to perform.
\begin{lemma}\label{lem:principe_maximum}$\left(\text{Maximum principe} \right)$.
  Under assumptions ({H}\ref{hyp:H1})-({H}\ref{hyp:H6}).  Let
  $(s_{\alpha,K}^{0})_{K \in \mathcal{T}}\in [0,1]$ and let $\D =
  \Big(\T,\E,(x_K)_{K\in\T},N,(t^n)_{n\in[0,N]}\Big)$ be a
  discretization of $\Omega\times (0,T)$ in the sense of Definition
  \ref{def:disc_QT} and assume that $(p_{\alpha,\D})$ is a solution of
  the finite volume \eqref{eq:p0}-\eqref{def_disc:pc}. Then, the
  saturation $(s_{\alpha,K}^{n})_{K \in \mathcal{T},n \in
    \{0,\ldots,N\}}$ remains in $[0,1]$.
\end{lemma}
\begin{proof}
  Let us show by induction in $n$ that for all $K \in \mathcal{T}, ~
  s^n_{\alpha,K}\geq 0$ where $\alpha=\2,\1$. For $\alpha=\2$, the
  claim is true for $n=0$ and for all $K \in \mathcal{T}$. We argue by
  induction that for all $K \in \mathcal{T}$, the claim is true up to
  order $n$. We consider the control volume $K$ such that
  $s^{n+1}_{\2,K}=\min{\{s^{n+1}_{\2,L}\}}_{L\in
    \mathcal{T}}$ and we seek that $s^{n+1}_{\2,K}\geq 0$.\\
  For the above mentioned purpose, multiply the equation in
  \eqref{sys_disc:pl} by $-(s_{\2,K}^{n+1})^-$, we obtain
  \begin{multline}\label{non-negative}
    - \abs{K}\phi_K\frac{\rho_\2(p^{n+1}_{\2,K})
      s^{n+1}_{\2,K}-\rho_\2(p^{n}_{\2,K}) s^{n}_{\2,K}}{\h}
    (s_{\2,K}^{n+1})^- \\ - \sum_{L \in N(K) } \tokl
    \rho^{n+1}_{\2,\KL}
    G_\2(s^{n+1}_{\2,K},s^{n+1}_{\2,L};\de(p_{\2})) (s_{\2,K}^{n+1})^-
    - F^{(n+1)}_{\2,K}(s_{\2,K}^{n+1})^- \\ - \abs{K}
    \rho_\2(p_{\2,K}^{n+1}) s_{\2,K}^{n+1}
    f_{P,K}^{n+1}(s_{\2,K}^{n+1})^- = -\abs{K} \rho_\2(p_{\2,K}^{n+1})
    \sgi f_{I,K}^{n+1}(s_{\2,K}^{n+1})^-\le 0.
\end{multline}
The numerical flux $G_\2$ is nonincreasing with respect to
$s_{\2,L}^{n+1}$ (see (a) in \eqref{Hypfluxes}), and consistence (see
(c) in \eqref{Hypfluxes}), we get

\begin{align}\label{nonnegat:I:2}
  G_\2(s^{n+1}_{\2,K},s^{n+1}_{\2,L};\de(p_{\2}))\,(s_{\2,K}^{n+1})^-
  &\le G_\2(s^{n+1}_{\2,K},s^{n+1}_{\2,K}
  ;\de(p_{\2}))\,(s_{\2,K}^{n+1})^-\notag\\
  &=-\de(p_{\2})\,M_\2(s^{n+1}_{\2,K}) \,(s_{\2,K}^{n+1})^-=0.
\end{align}
Using the identity
$s_{\2,K}^{n+1}=({s_{\2,K}^{n+1}})^+-(s_{\2,K}^{n+1})^-$,
and the mobility $M_\2$ extended by zero on $]-\infty, 0]$, then
$M_\2(s^{n+1}_{\2,K}) (s_{\2,K}^{n+1})^- = 0$ and
\begin{multline}\label{nonnegat:I:3}
  -F^{(n+1)}_{\2,K}(s_{\2,K}^{n+1})^- -\abs{K} \rho_\2(p_{\2,K}^{n+1})
  s_{\2,K}^{n+1} f_{P,K}^{n+1}(s_{\2,K}^{n+1})^- \\ = \sum_{L \in N(K)
  } (\rho_{\2,\KL}^{n+1})^2 M_\2(s^{n+1}_{\2,L}){\bf
    g}_{L,K}(s_{\2,K}^{n+1})^- +\abs{K} \rho_\2(p_{\2,K}^{n+1})
  f_{P,K}^{n+1}((s_{\2,K}^{n+1})^-)^2 \geq 0.
\end{multline}
Then, we deduce from \eqref{non-negative} that
$$  
\rho_\2(p^{n+1}_{\2,K})|(s^{n+1}_{\2,K})^{-}|^2
  +\rho_\2(p^{n}_{\2,K})s^{n}_{\2,K}(s_{\2,K}^{n+1})^-\leq
0,
$$
and from the nonnegativity of $s^{n}_{\2,K}$, we obtain
$(s_{\2,K}^{n+1})^-=0$.  This implies that $s_{\2,K}^{n+1}\geq
0$ and
$$
0\le s^{n+1}_{\2,K} \le s^{n+1}_{\2,L} \text{ for all } n \in [0,N-1]
\text{ and }L \in \mathcal{T}.
$$ 
In the same way, we prove $s_{\1,K}^{n+1}\ge 0$.
\end{proof}
\subsection{Estimations on the pressures}
\begin{prop}\label{prop:estimation_pression}
  Let $p_{\alpha,\D}$ be
  a solution of \eqref{eq:p0}-\eqref{def_disc:pc}. Then, there exists
  a constant $C>0$, which only depends on $M_\alpha$, $\Omega$, $T$,
  $p^{0}_\alpha$, $s^0_\alpha$, $s^I_\alpha$, $f_P$, $f_I$ and not on
  $\D$, such that the following discrete $L^2(0,T;H^1(\Omega))$
  estimates hold:
\begin{align}\label{est:p_alpha}
  \sm\tokl M_\alpha(s_{\alpha,K|L}^{n+1})
  |p_{\alpha,L}^{n+1}-p_{\alpha,K}^{n+1}|^2\le C,
\end{align}
and
\begin{align}\label{est:p_\2lobale}
  \sm \tokl|p_L^{n+1}-p_K^{n+1}|^2\le C.
\end{align}
\end{prop}
\begin{proof} We define the function $
  \mathcal{H}_{\alpha}(p_{\alpha}) :=
  \rho_{\alpha}(p_{\alpha})g_{\alpha}(p_{\alpha}) - p_{\alpha},$ $
  \pc(s_\3) := \int_{0}^{s_\3}p_c(z)\dd z $ and $
  g_{\alpha}(p_{\alpha}) =
  \int_{0}^{p_{\alpha}}\frac{1}{\rho_{\alpha}(z)}\dd z$. In the
  following proof, we denote by $C_i$ various real values which only
  depend on $M_\alpha$, $\Omega$, $T$, $p^{0}_\alpha$, $s^0_\alpha$,
  $s^I_\alpha$, $f_P$, $f_I$ and not on $\D$. To prove the estimate
  \eqref{est:p_alpha}, we multiply (\ref{sys_disc:pl}) and
  (\ref{sys_disc:pg}) respectively by $g_\2(p_{\2,K})$,
  $g_\1(p_{\1,K})$ and adding them, then summing the resulting
  equation over $K$ and $n$.  We thus get:
\begin{equation}\label{disc_estimation}
  E_{1}+E_{2}+E_{3}+E_4= 0,
  \end{equation}
where
\begin{align*}
  &E_{1} = \sn \sum_{K \in \mathcal{T}} \abs{K}\phi_K \Big(
  (\rho_\2(p^{n+1}_{\2,K})s^{n+1}_{\2,K}
  -\rho_\2(p^{n}_{\2,K})s^{n}_{\2,K})\; g_\2(p^{n+1}_{\2,K}) \\ &
  \hspace{8.3cm} + (\rho_\1(p^{n+1}_{\1,K})s^{n+1}_{\1,K}
  -\rho_\1(p^{n}_{\1,K})s^{n}_{\1,K})\; g_\1(p^{n+1}_{\1,K})
  \Big),\\
  &E_{2} = \sn\pas \sum_{K \in \mathcal{T}} \sum_{L \in N(K)}\tokl
  \Big(\rho^{n+1}_{\2,\KL}
  G_\2(s^{n+1}_{\2,K},s^{n+1}_{\2,L};\de(p_{\2}))\;
  g_\2(p_{\2,K}^{n+1})\\
  &\hspace{9.3cm} + G_\1(s^{n+1}_{\1,K},s^{n+1}_{\1,L};\de(p_{\1}))\;
  g_\1(p_{\1,K}^{n+1})
  \Big),\\
  & E_{3}=\sn\pas \sum_{K \in \mathcal{T}} \Big( F^{(n+1)}_{\2,\KL}\;
  g_\2(p_{\2,K}^{n+1})+ F^{(n+1)}_{\1,\KL}\; g_\1(p_{\1,K}^{n+1})
  \Big),\\
  & E_4 = \sn\pas \sum_{K \in \mathcal{T}}\abs{K}
  \Big(\rho_\2(p_{\2,K}^{n+1}) s_{\2,K}^{n+1} f_{P,K}^{n+1}
  g_\2(p_{\2,K}^{n+1}) - \rho_\2(p_{\2,K}^{n+1}) \sgi f_{I,K}^{n+1}
  g_\2(p_{\2,K}^{n+1}) \\ & \hspace{3.5cm} + \rho_\1(p_{\1,K}^{n+1})
  s_{\1,K}^{n+1} f_{P,K}^{n+1} g_\1(p_{\1,K}^{n+1}) -
  \rho_\1(p_{\1,K}^{n+1}) \sgi f_{I,K}^{n+1} g_\1(p_{\1,K}^{n+1})
  \Big).
\end{align*}
To handle the first term of the equality \eqref{disc_estimation}. Let us 
forget the exponent $n+1$ and let note with the exponent $*$ the
physical quantities at time $t^n$.  In \cite{ZS10} the authors prove that  : 
for all $s_\alpha\ge 0$ and $s^\star_\alpha\ge 0$ such that
$s_\2+s_\1=s^\star_\2 +s^\star_\1=1$,
\begin{multline}
  \label{bibi12}
  \bigl(\rho_\2(p_\2)s_\2-\rho_\2(p_\2^\star)
  s^\star_\2\bigr)g_\2(p_\2)+
  \bigl(\rho_\1(p_\1)s_\1-\rho_\1(p_\1^\star)
  s^\star_\1\bigr)g_\1(p_\1) \\
  \ge\mathcal{H}_\2(p_\2)s_\2-\mathcal{H}_\2(p_\2^\star)s^\star_\2 +
  \mathcal{H}_\1(p_\1)s_\1-\mathcal{H}_\1(p_\1^\star)s^\star_\1 -
  \pc(s_\3)+\pc(s_\3^\star).
\end{multline}
The proof of  \eqref{bibi12} is based on the concavity property
of $g_\alpha$ and    $\pc$.
 So, this yields to 
\begin{multline}\label{E1}
  E_1\geq  \sum_{K\in \mathcal{T}} 
      \phi_K \abs{K} \Big(s_{\2,K}^N \mathcal{H}(p_{\2,K}^N) 
                        -  s_{\2,K}^0 \mathcal{H}(p_{\2,K}^0)
                        +  s_{\1,K}^N \mathcal{H}(p_{\1,K}^N) 
                        -  s_{\1,K}^0 \mathcal{H}(p_{\1,K}^0)
                    \Big) \\  
 - \sum_{K\in \mathcal{T}}\phi_K \abs{K} \pc(s_{\3,K}^N) 
 + \sum_{K\in \mathcal{T}} \phi_K \abs{K}\pc(s_{\3,K}^0).
\end{multline}
Using the fact that the numerical fluxes $G_\2$ and $G_\1$ are
conservative in the sense of (c) in (\ref{Hypfluxes}), we obtain by
discrete integration by parts (see Lemma \ref{lem:integration_parti})
\begin{align*}
  E_2 = \frac{1}{2}\sn\pas \sum_{K \in \mathcal{T}} \sum_{L \in N(K)}\tokl 
    \Big(& \rho_{\2,\KL}^{n+1} G_\2(s^{n+1}_{\2,K},s^{n+1}_{\2,L};\de(p_{\2}))            (g_\2(p_{\2,K}^{n+1})-g_\2(p_{\2,L}^{n+1})) \\ &
        +  \rho_{\1,\KL}^{n+1} G_\1(s^{n+1}_{\1,K},s^{n+1}_{\1,L};\de(p_{\1}))            (g_\1(p_{\1,K}^{n+1})-g_\1(p_{\1,L}^{n+1}))
    \Big),
 \end{align*}
 and due to the correct choice of the density of the phase $\alpha$ on each interface,
\begin{align}\label{choice_of_density}
  \rho_{\alpha,\KL}^{n+1} (g_{\alpha}(p_{\alpha,K}^{n+1})-g_{\alpha}(p_{\alpha,L}^{n+1}))=
p_{\alpha,K}^{n+1}-p_{\alpha,L}^{n+1},
\end{align}
we obtain 
\begin{align*}
  E_2 = \frac{1}{2}\sn\pas \sum_{K \in \mathcal{T}} \sum_{L \in N(K)}\tokl
       & \Big( G_\2(s^{n+1}_{\2,K},s^{n+1}_{\2,L};\de(p_{\2}))
               (p_{\2,K}^{n+1}-p_{\2,L}^{n+1})\\ & +
              G_\1(s^{n+1}_{\1,K},s^{n+1}_{\1,L};\de(p_{\1}))
               (p_{\1,K}^{n+1}-p_{\1,L}^{n+1}) \Big).
\end{align*}
The definition of the upwind fluxes in \eqref{Gupwind} implies 
\begin{multline*}  
  G_\2(s^{n+1}_{\2,K},s^{n+1}_{\2,L};\de(p_{\2}))
        (p_{\2,K}^{n+1}-p_{\2,L}^{n+1})
     +
  G_\1(s^{n+1}_{\1,K},s^{n+1}_{\1,L};\de(p_{\1}))
        (p_{\1,K}^{n+1}-p_{\1,L}^{n+1})
\\ =   M_{\2}(s_{\2,\KL}^{n+1})(\de(p_{\2}))^2 
     +
       M_{\1}(s_{\1,\KL}^{n+1})(\de(p_{\1}))^2.
\end{multline*}
Then, we obtain the following equality
\begin{align}\label{E2}
  E_2 = \frac{1}{2}\sn\pas \sum_{K \in \mathcal{T}} \sum_{L \in
    N(K)} \tokl \Big(M_{\2}(s_{\2,\KL}^{n+1})(\de(p_\2))^2 +
  M_{\1}(s_{\1,\KL}^{n+1})(\de(p_\1)^2\Big).
\end{align}
To handle the other terms of the equality (\ref{disc_estimation}),
firstly let us remark that the numerical fluxes of gravity term are  conservative which satisfy 
$F_{\2,\KL}^{n+1}= -F_{\2,L,K}^{n+1}$ and $F_{\1,\KL}^{n+1}=
-F_{\1,L,K}^{n+1}$, so we integrate by parts and we obtain
$$
E_{3}= \frac1 2\sn\pas \sum_{K \in \mathcal{T}}\sum_{L \in N(K)}
|\sigma_{\KL}| \Big ( F^{(n+1)}_{\2,\KL} (g_\2(p_{\2,K}^{n+1}) -
                                          g_\2(p_{\2,L}^{n+1})) 
                    + F^{(n+1)}_{\1,\KL} (g_\1(p_{\1,K}^{n+1}) -
                                          g_\1(p_{\1,L}^{n+1})) 
               \Big).
$$
According to the choice of the density of the phase $\alpha$ on each
interface \eqref{choice_of_density} and the definition \eqref{fgravity} we obtain
\begin{align*}
  E_3 = & - 
   \frac1 2\sn\pas \sum_{K \in \mathcal{T}}\sum_{L \in
     N(K)}\abs{\sigma_{\KL}}\rho^{n+1}_{\2,\KL}[M_\2(s_{\2,K}^{n+1})
     \G^+_{\KL}-M_\2(s_{\2,L}^{n+1}) \G^-_{\KL}](\de(p_\2))
      \\ &\quad - 
   \frac1 2\sn\pas \sum_{K \in \mathcal{T}}\sum_{L \in
     N(K)}\abs{\sigma_{\KL}}\rho^{n+1}_{\1,\KL}[M_\1(s_{\1,K}^{n+1})
     \G^+_{\KL}-M_\1(s_{\1,L}^{n+1}) \G^-_{\KL}](\de(p_\1)).
\end{align*}


Recall the truncations of $\de(p_\alpha)$ 
$$
(\de(p_\alpha))^+ = max\{\de(p_\alpha),0\},\;\quad  (\de(p_\alpha))^-
= max\{- \de(p_\alpha),0\},
$$
with $\de(p_\alpha) = (\de(p_\alpha))^+ - (\de(p_\alpha))^-.$ So we obtain
\begin{align*}
  E_3 &\le \frac1 2\sn\pas \sum_{K \in \mathcal{T}}\sum_{L \in N(K)} 
            \abs{\sigma_{\KL}} \rho^{n+1}_{\2,\KL} M_\2(s_{\2,K}^{n+1}) 
                                \G^+_{\KL}(\de(p_\2))^-\\ &\quad 
         + \frac1 2\sn\pas \sum_{K \in \mathcal{T}}\sum_{L \in N(K)}
            \abs{\sigma_{\KL}} \rho^{n+1}_{\2,\KL} M_\2(s_{\2,L}^{n+1}) 
                                \G^-_{\KL}(\de( p_\2))^+ \\ & \quad
         + \frac1 2\sn\pas \sum_{K \in \mathcal{T}}\sum_{L \in N(K)}
            \abs{\sigma_{\KL}}\rho^{n+1}_{\1,\KL}M_\1(s_{\1,K}^{n+1})
                                \G^+_{\KL}(\de(p_\1))^- \\ &\quad 
         + \frac1 2\sn\pas \sum_{K \in \mathcal{T}}\sum_{L \in N(K)}
            \abs{\sigma_{\KL}}\rho^{n+1}_{\1,\KL}M_\1(s_{\1,L}^{n+1})
                                \G^-_{\KL}(\de(p_\1))^+. 
\end{align*}
From the following equality 
$\abs{\sigma_{\KL}}=(\dis |\sigma_{\KL}|)^{\frac 1
  2}\tokl^{\frac 1 2}$ and
apply the Cauchy-Schwarz inequality to obtain
\begin{align*}
  E_3 \le & 2C \sn\pas \sum_{K \in \mathcal{T}}\sum_{L \in N(K)}
  \dis  |\sigma_{\KL}| \\ & + \frac 1 4 \sn\pas \sum_{K \in
    \mathcal{T}}\sum_{L \in N(K)} \tokl\Big(
  M_{\2}(s_{\2,L}^{n+1})\; ((\de(p_\2))^+)^2 + M_{\2}(s_{\2,K}^{n+1})\;
  ((\de(p_\2))^-)^2 \\ & \qquad\qquad\qquad\qquad\qquad\qquad\quad +
  M_{\1}(s_{\1,L}^{n+1})\; ((\de(p_\1))^+)^2 + M_{\1}(s_{\1,K}^{n+1})\;
  ((\de(p_\1))^-)^2\Big)\\ & \le 2CT|\Omega| \\ & + \frac 1 4
  \sn\pas \sum_{K \in \mathcal{T}}\sum_{L \in N(K)}
  \tokl \Big(M_{\2}(s_{\2,L}^{n+1})\;
  ((\de(p_\2))^+)^2 + M_{\1}(s_{\2,K}^{n+1})\; ((\de(p_\2))^-)^2 \\ &
  \qquad\qquad\qquad\qquad\qquad\qquad\quad+ M_{\1}(s_{\1,L}^{n+1})\;
  ((\de(p_\1))^+)^2 + M_{\1}(s_{\1,K}^{n+1})\; ((\de(p_\1))^-)^2\Big).
\end{align*}
From the definition of the truncations of $\de(p_\alpha)$, we obtain
\begin{multline}\label{E3}
  E_3 \le 2CT|\Omega|+ \frac 1 4 \sn\pas \sum_{K \in
    \mathcal{T}}\sum_{L \in N(K)}
  \tokl\Big(
  M_{\2}(s_{\2,\KL}^{n+1})(\de(p_\2))^2 \\ +
  M_{\1}(s_{\1,\KL}^{n+1})(\de(p_\1))^2\Big).
\end{multline}
The last term will be absorbed by the terms on pressures from the
estimate \eqref{E2}.\\
In order to estimate $E_4$, using the fact that the densities are
bounded and the map $g_\alpha$ is sublinear $(\text{a.e.}|g(p_\alpha)|\le
C |p_\alpha|)$, we have
\begin{equation*}
  \abs{E_4} \le C_1\sn\pas \sum_{K \in
    \mathcal{T}}\abs{K}(f_{P,K}^{n+1}+f_{I,K}^{n+1})
  (|p_{\2,K}^{n+1}|+|p_{\1,K}^{n+1}|),
\end{equation*}
then
\begin{equation*}
  \abs{E_4} \le C_1\sn\pas \sum_{K \in
    \mathcal{T}}\abs{K}(f_{P,K}^{n+1}+f_{I,K}^{n+1})
  (2|p_{K}^{n+1}|+|\bar{p}_{K}^{n+1}|+|\tilde{p}_{K}^{n+1}|).
\end{equation*}
Hence, by the H\"older inequality, we get that
\begin{equation*}
  \abs{E_4} \leq C_2\norm{f_P+f_I}_{L^2(Q_T)} \big(
  \sum_{n=0}^{N-1}\pas
  \norm{p^{n+1}_h}_{L^2({\Omega})}^2\big)^{\frac{1}{2}},
\end{equation*}
and, from the discrete Poincar\'e inequality lemma \ref{lem:poincare},
we get
\begin{equation}\label{E4}
  \abs{E_{4}} \leq C_3 \big( \sum_{n=0}^{N-1}\pas
  \norm{p^{n+1}_h}_{H_h}^2\big)^{\frac{1}{2}}+C_4.
\end{equation}
The equality \eqref{disc_estimation} with the inequalities \eqref{E1},
\eqref{E2}, \eqref{E3}, \eqref{E4} give \eqref{est:p_alpha}. Then we
deduce \eqref{est:p_\2lobale} from \eqref{ineq:pression_globale}.
\end{proof}

We now state the following corollary, which is essential for the
compactness and limit study.
\begin{cor}\label{cor:est} From the previous Proposition, we deduce
  the following estimations:
  \begin{align}
    &\sm \tokl (\de(\B(s_\3)))^2 \le C,\label{est:discrete_beta}\\ &
    \sm \tokl M_{\2,K|L}^{n+1}(\de(\bar{p}(s_\3)))^2 \le C,
    \label{est:discrete_pbar}
    \end{align}
    and
\begin{align}
  \sm \tokl M_{\1,K|L}^{n+1}(\de(\tilde{p}(s_\3)))^2 \le
  C.\label{est:discrete_ptilde}
  \end{align}
\end{cor}
\begin{proof}
  The prove of the estimates \eqref{est:discrete_beta},
  \eqref{est:discrete_pbar} and \eqref{est:discrete_ptilde} are a
  direct consequence of the inequality \eqref{ineq:discrete_beta},
  \eqref{ineq:discrete_pbar}, \eqref{ineq:discrete_ptilde} and the
  Proposition \ref{prop:estimation_pression}.
\end{proof}

\section{Existence of the finite volume scheme}\label{sec:existence}
We start with a technical assertion to characterize the zeros of a vector field which stated and proved in \cite{evans:book}.  
\begin{lemma}\label{lem:exist-classic}(\cite{evans:book}, p. 529)
  Assume the continuous function $v : \R^n \rightarrow \R^n$ satisfies 
  $$
  v(z)\cdot z \ge 0 \text{ if } \|z\| = r,
  $$
  for some $r>0$. Then there exists a point $z$ with $\|z\|\le r$ such that 
  $$
  v(z)=0.
  $$
\end{lemma}

\begin{prop}
\label{prop:existance}
The problem \eqref{sys_disc:pl}-\eqref{sys_disc:pg} admits at least one solution
$(p^{n}_{\2,K},p^{n}_{\1,K})_{(K,n) \in \D}$.
\end{prop}
\begin{proof}
  At the beginning of the proof, we set the following notations;
\begin{align*}
  &        \mathcal{M}:=Card(\mathcal{T}),\\ &
    p_{\2,\mathcal{M}}:=\{p^{n+1}_{\2,K} \}_{K\in\mathcal{T}}
                                         \in\R^{\mathcal{M}},\\ & 
    p_{\1,\mathcal{M}}:=\{p^{n+1}_{\1,K} \}_{K\in\mathcal{T}} 
                                         \in\R^{\mathcal{M}}.
\end{align*}
We define the map
$\mathcal{T}_h: \R^\mathcal{M} \times \R^\mathcal{M} 
                  \longrightarrow
                \R^\mathcal{M} \times \R^\mathcal{M},$

$$\mathcal{T}_h(p_{\2,\mathcal{M}}, p_{\1,\mathcal{M}})=
          (\{\mathcal{T}_{\2,K}\}_{K\in \mathcal{T}},
            \{\mathcal{T}_{\1,K}\}_{K\in \mathcal{T}})
\,\, \text{where,} $$
\begin{align}
  & \mathcal{T}_{\2,K} = \abs{K}\phi_K 
                          \frac{\rho_\2(p^{n+1}_{\2,K})s^{n+1}_{\2,K}-
                                \rho_\2(p^{n}_{\2,K})s^{n}_{\2,K}}{\h} 
                       + \sum_{L \in N(K) } \tokl \rho^{n+1}_{\2,\KL}
                         G_\2(s^{n+1}_{\2,K},s^{n+1}_{\2,L};\de(p_\2))
  \notag \\ & \hspace{5.5cm}
                       +F^{\;n+1}_{\2,K} 
                       +\abs{K}\rho_\2(p_{\2,K}^{n+1})
                         \big( 
                              s_{\2,K}^{n+1} f_{P,K}^{n+1} -
                              \sli f_{I,K}^{n+1} 
                         \big),
  \\ & 
   \mathcal{T}_{\1,K} = \abs{K}\phi_K 
                          \frac{\rho_\1(p^{n+1}_{\1,K})s^{n+1}_{\1,K}-
                                \rho_\1(p^{n}_{\1,K})s^{n}_{\1,K}}{\h} 
                       + \sum_{L \in N(K) }\tokl \rho^{n+1}_{\1,\KL}
                         G_\1(s^{n+1}_{\1,K},s^{n+1}_{\1,L};\de(p_\1))
  \notag \\ & \hspace{5.5cm}
                       +F^{\;n+1}_{\1,K} 
                       +\abs{K}\rho_\1(p_{\1,K}^{n+1})
                         \big( 
                              s_{\1,K}^{n+1} f_{P,K}^{n+1} -
                              \sgi f_{I,K}^{n+1} 
                         \big).
\end{align}
Note that $\mathcal{T}_h$ is well defined as a continuous function.
Also we define the following homeomorphism
$\mathcal{F}:\R^\mathcal{M} \times \R^\mathcal{M} \mapsto
             \R^\mathcal{M} \times \R^\mathcal{M}$ 
such that,
$$
\mathcal{F}(p_{\2,\mathcal{M}},p_{\1,\mathcal{M}})=
 (v_{\2,\mathcal{M}}, v_{\1,\mathcal{M}})
$$
where $v_{\alpha,\mathcal{M}}=
      \{g_\alpha(p^{n+1}_{\alpha,K})\}_{K\in \mathcal{T}}.$\\
Now let us consider the following continuous mapping $\mathcal{P}_h$
defined as
\begin{align*}
  \mathcal{P}_h(v_{\2,\mathcal{M}}, v_{\1,\mathcal{M}})& = 
  \mathcal{T}_h \circ \mathcal{F}^{-1}(v_{\2,\mathcal{M}},v_{\1,\mathcal{M}})=
  \mathcal{T}_h(p_{\2,\mathcal{M}},p_{\1,\mathcal{M}}).
\end{align*}
According to Lemma \ref{lem:exist-classic}, our goal now is to show that 
\begin{align} \mathcal{P}_h(v_{\2,\mathcal{M}},v_{\1,\mathcal{M}}) \cdot 
                            (v_{\2,\mathcal{M}},v_{\1,\mathcal{M}}) > 0,
            \quad \text{ for }
               \norm{(v_{\2,\mathcal{M}},v_{\1,\mathcal{M}})}_{\R^{2\mathcal{M}}}=r>0,\label{goal-exis}
\end{align}
and for a sufficiently large $r$.\\
We observe that
\begin{equation*}
\begin{split}
    \mathcal{P}_h(v_{\2,\mathcal{M}},v_{\1,\mathcal{M}}) \cdot 
                  (v_{\2,\mathcal{M}},v_{\1,\mathcal{M}})   \ge &
     \frac{1}{\h}\sum_{K\in \mathcal{T}} \phi_K \abs{K}
     \Big(
       s_{\2,K}^{n+1} \mathcal{H}(p_{\2,K}^{n+1}) - 
       s_{\2,K}^n     \mathcal{H}(p_{\2,K}^n)       \\ & \hspace{3cm} +
       s_{\1,K}^{n+1} \mathcal{H}(p_{\1,K}^{n+1}) - 
       s_{\1,K}^n     \mathcal{H}(p_{\1,K}^n)
     \Big) \\ & -\frac{1}{\h}\pc(s_{\3,K}^{n+1}) +
  \frac{1}{\h}\pc(s_{\3,K}^n) +C \norm{p^{n+1}_{h}}^2_{H_h(\Omega)}
  - C,   
\end{split}
\end{equation*}
for some constants $C>0$. This implies that
\begin{equation}\label{eq:disc-exist:5}
  \begin{split}
    \mathcal{P}_h(v_{\2,\mathcal{M}},v_{\1,\mathcal{M}}) \cdot 
                  (v_{\2,\mathcal{M}},v_{\1,\mathcal{M}})  \ge & 
   - \frac{1}{\pas} \sum_{K\in \mathcal{T}} \phi_K \abs{K}
      \Big(
           s_{\2,K}^n \mathcal{H}(p_{\2,K}^n) +
           s_{\1,K}^n \mathcal{H}(p_{\1,K}^n) 
      \Big)   \\  & 
   -\frac{1}{\pas}\pc(s_{\3,K}^{n+1})+C\norm{p^{n+1}_{h}}^2_{H_h(\Omega)} -C',
\end{split}
\end{equation} 
for some constants $C,C'>0$.  Finally using the fact that $g_\alpha$ is a
Lipschitz function, then there exists a constant $C>0$ such that
\begin{align*}
  \norm{(\{g_\2(p^{n+1}_{\2,K})\}_{K\in \mathcal{T}},
         \{g_\1(p^{n+1}_{\1,K}))\}_{K\in \mathcal{T}})}_{\R^{2\mathcal{M}}}&\le
           C \Big(\norm{p^{n+1}_{\2,h}}_{L^2(\Omega)} +
     \norm{p^{n+1}_{\1,h}}_{L^2(\Omega)}\Big) \\ & \le
      2C
  \Big(\norm{p^{n+1}_{h}}_{L^2(\Omega)}
  +\norm{\bar{p}^{n+1}_{h}}_{L^2(\Omega)} +
  \norm{\tilde{p}^{n+1}_{h}}_{L^2(\Omega)}\Big) \\ & \le 2C
  \big(\norm{p^{n+1}_{h}}_{H_h(\Omega)} +C_1\big).
\end{align*}
Using this to deduce from \eqref{eq:disc-exist:5} that
(\ref{goal-exis}) holds for $r$ large enough.  Hence, we obtain the
existence of at least one solution to the scheme
\eqref{sys_disc:pl}-\eqref{sys_disc:pg}.
\end{proof}

\section{Compactness properties}\label{sec:compacity}
In this section we derive estimates on differences of space and time
translates of the function $\phi_{\D} \rho_\alpha(p_{\alpha,\D}) s_{\alpha,\D}$
which imply that the sequence $\phi_{\D} \rho_\alpha(p_{\alpha,\D})s_{\alpha,\D}$
is relatively compact in $L^1(Q_T)$.

We replace the study of discrete functions $U_{\alpha,\D}=\phi_{\D}
\rho_\alpha(p_{\alpha,\D}) s_{\alpha,\D}$ (constant per cylinder
$\QKn:=(t^n,t^{n+1})\times K$) by the study of functions $\bar
U_{\alpha,\D}=\phi_\D \rho_\alpha(\bar p_{\alpha,\D}) \bar
s_{\alpha,\D}$ piecewise continuous in $t$ for all $x$, constant in
$x$ for all volume $K$, defined as
$$
\bar U_{\alpha,\D}(t,x)=\snk \frac
1{\delta t}\Bigl(\;(t-n\delta t)U^{n+1}_{\alpha,K} \,+\,((n+1)\delta
t-t)U^{n}_{\alpha,K} \; \Bigr)\;\charu_{Q_K^n}(t,x).
$$

One may deduce from the estimates \eqref{est:p_\2lobale} and
\eqref{est:discrete_beta} the following property.
\begin{lemma}\label{lem:translater-espace}$\left(\text{Space translate
      of } \bar U_{\alpha,\D}\right)$. Under the assumptions
  $({H}\ref{hyp:H1})-({H}\ref{hyp:H6})$ . Let $\D$ be a finite volume
  discretization of $\Omega\times (0,T)$ in the sense of Definition
  \ref{def:disc_QT} and let $p_{\alpha,\D}$ be a solution of
  \eqref{eq:p0}-\eqref{def_disc:pc}. Then, the following inequality
  hold:
  \begin{equation}\label{eq:translater_espace}
    \int_{\Omega^{'}\times (0,T)}\abs{\bar U_{\alpha,\D}(t,x+y) - 
                                      \bar U_{\alpha,\D}(t,x)}  
                                    \dd x \dd t\le \omega(\abs{y}),
\end{equation}
for all $y \in \R^\dm$ with $\Omega '=\{x \in \Omega, \, [x,x+y]\subset \Omega\}$ and $\omega(\abs{y}) \to 0$ when $\abs{y}\to 0$.
\end{lemma}
\begin{proof}
For $\alpha = \2$ and from the definition of $U_{\2,\D}$, one gets

\begin{equation*}
      \begin{split}
        &\int_{(0,T)\times\Omega^{'}}\abs{U_{\2,\D}(t,x+y)-U_{\2,\D}(t,x)}
        \dd x\dd t\\
        & = \int_{(0,T)\times\Omega^{'}}
        \abs{\Big(\rho_\2(p_{\2,{\D}})s_{\2,\D}\Big)(t,x+y) -
          \Big(\rho_\2(p_{\2,{\D}})s_{\2,\D}\Big)(t,x)} \dd x\dd t\\
        & \le \int_{(0,T)\times\Omega^{'}}
        \abs{s_{\2,\D}(t,x+y) \Big(\rho_\2(p_{\2,{\D}}(t,x+y)) -\rho_\2(p_{\2,{\D}}(t,x))\Big) }\dd x\dd t \\
        & \quad +   \int_{(0,T)\times\Omega^{'}} \abs{ \rho_\2(p_{\2,{\D}})(t,x)\Big(s_{\2,\D}(t,x+y)  - s_{\2,\D}(t,x)\Big) }\dd x\dd t\\
        & \le E_1+E_2
        \end{split}
\end{equation*}
where $E_1$ and $E_2$ defined as follows
\begin{equation}\label{e1trans}
E_1 = \rho_M \int_{(0,T)\times\Omega^{'}} \abs{s_{\2,\D}(t,x+y)  - s_{\2,\D}(t,x) }\dd x\dd t,
\end{equation}
\begin{equation}\label{e2transs}
E_2 = \int_{(0,T)\times\Omega^{'}}
        \abs{ \rho_\2(p_{\2,{\D}}(t,x+y)) -\rho_\2(p_{\2,{\D}}(t,x)) }\dd x\dd t.
\end{equation}
To handle with the  space translation on saturation, we use the fact that $\B^{-1}$ is an h\"older function, then 
$$
E_1\le  \rho_M C \int_{(0,T)\times\Omega^{'}}
        \abs{\B(s_{\2,\D}(t,x+y) )- \B(s_{\2,\D}(t,x))}^\theta \dd x\dd t
$$
and by application of the Cauchy-Schwarz inequality, we deduce
$$
E_1\le  C \Big(\int_{(0,T)\times\Omega^{'}}
        \abs{\B(s_{\2,\D}(t,x+y) )- \B(s_{\2,\D}(t,x))} \dd x\dd t\Big)^\theta. 
$$

According to \cite{Eymard:book}), let $y \in \R^\dm$, $x \in \Omega '$, and $L \in N(K)$. We set
$$
\beta_{\sigma_{\KL}}=
\begin{cases}
  1, & \text{if the line segment $[x,x+y]$ intersects $\sigma_{\KL}$, 
    $K$ and $L$},\\
  0, & \text{otherwise}.
\end{cases}
$$
We observe that (see for more details \cite{Eymard:book})
\begin{equation}\label{est:space}\begin{split}
    &\int_{\Omega '}\beta_{\sigma_{\KL}}(x) \,dx \le |\sigma_{\KL}|
    \abs{y}.
\end{split}\end{equation}
To simplify the notation, we write $\underset{\sigma_{\KL}}\sum$ instead of $\underset{\{(K,L)\in
  {\mathcal{T}}^2,\,K\ne L,\,\abs{\sigma_{\KL}}\ne 0\}}\sum$.

Now, denote that
\begin{equation*}
      \begin{split}
        E_1 \le & C \Big(\sn \h \sum_{\sigma_{K,L}} 
              \Big|\B(s_{\2,L}) - \B(s_{\2,K})\Big|  
               \int_{\Omega '}\beta_{\sigma_{\KL}}(x)\dd x\Big)^\theta \\
        & \le 
            C \Big(\abs{y} \sn \h \sum_{\sigma_{K,L}} \abs{\sigma_{\KL}}
              \Big|\B(s_{\2,L}) - \B(s_{\2,K})\Big|\Big)^\theta .
\end{split}
\end{equation*}
Let us again write  $\abs{\sigma_{K,L}} = (d_{K,L}|\sigma_{K,L}|)^{\frac 1 2}\tokl^{\frac 1 2}$, applying again the Cauchy-Schwarz inequality  and using the fact that the discrete gradient of the function $\B$ is bounded \eqref{est:discrete_beta} to obtain
\begin{equation} \label{e1trans2}
 E_1 \le C \abs{y}^\theta. 
\end{equation}

To treat the space translate of $E_2$, we use the fact that the map
$\rho_\2^\prime$ is bounded and the relationship between the gas
pressure and the global pressure, namely : $p_\2=p-\bar{p}$ defined
in \eqref{def:pression_globale}, then we have
\begin{equation}\label{e2trans}
   \begin{split}
E_2 &\le  \max_\R |\rho_\2^\prime |\int_{(0,T)\times\Omega^{'}}
       \abs{p_{\2,{\D}}(t,x+y) -p_{\2,{\D}}(t,x)}\dd x\dd t\\
         & \le  
          \max_\R |\rho_\2^\prime  | \int_{(0,T)\times\Omega^{'}}
       \abs{p_{\D}(t,x+y) -p_{\D}(t,x)}\dd x\dd t \\
       &+ 
         \max_\R |\rho_\2^\prime  | \int_{(0,T)\times\Omega^{'}}\abs{\bar{p}(s_{\2,\D}(t,x+y)) -\bar{p}(s_{\2,\D}(t,x)) }\dd x\dd t,
        \end{split}
\end{equation}
furthermore one can easily show that $\bar{p}$ is a $C^1([0,1];\R)$,
it follows, there exists a positive constant $C>0$ such that
\begin{equation*}
      \begin{split}
        E_2&\le C \int_{(0,T)\times\Omega^{'}}
            | p_{\D}(t,x+y) - p_{\D}(t,x) | \dd x \dd t\\ 
        & \qquad + C \int_{(0,T)\times\Omega^{'}}
        |s_{\2,\D}(t,x+y) -  s_{\2,\D}(t,x)| \dd x \dd t.
\end{split}
\end{equation*}
The last term in the previous inequality is proportional to $E_1$, and
consequently it remains to show that the space translate on the global
pressure is small with $y$. In fact
\begin{equation*}
      \begin{split}
           \int_{(0,T)\times\Omega^{'}}  |p_{\D}(t,x+y) - p_{\D}(t,x) | \dd x \dd t 
               & \le \sn \h \sum_{\sigma_{K,L}} |p_{L}^{n+1} - p_{K}^{n+1} | 
            \int_{\Omega '}\beta_{\sigma_{\KL}}(x)\dd x \\
            & \le \abs{y} \sn \h \sum_{\sigma_{K,L}} \abs{\sigma_{\KL}}
               |p_{L}^{n+1} - p_{K}^{n+1}|.
\end{split}
\end{equation*}
Finally, using the fact that the discrete gradient of global pressure
is bounded \eqref{est:p_\2lobale}, we deduce that
 \begin{equation}\label{est:space:4}
      \begin{split}
        &\int_{(0,T)\times\Omega^{'}}\abs{U_{\2,\D}(t,x+y)-U_{\2,\D}(t,x)}
        \dd x \le C (\abs{y} +\abs{y}^\theta),
\end{split}
\end{equation}
for some constant $C>0$.\\
In addition, we have
\begin{align*} \int_0^{+\infty}\int_{\Omega^{'}}
  |\bar{U}_{\2,\D}(t,x+\dd x)- \bar{U}_{\2,\D}(t,x) |\dd x \dd t \leq &2
  \int_0^{T}\int_{\Omega^{'}}
  |U_{\2,\D}(t,x+\dd x) - U_{\2,\D}(t,x) |\dd x \dd t \\
  +&2\delta t\int_{\Omega^{'}_\delta}|U_{\2,\D}^0(x)|\,\dd x
\end{align*}
where $U_\2^0=\rho_\2(p_{\2}^0)s_{\2}^0$ and $\Omega^{'}_\delta=\{x
\in \Omega, \, \text{dist}(x,\Omega^{'})< \abs{\delta}\}$. By
\eqref{est:space:4}, the assumption $\h \to 0$ as $\size(\D)\to 0$ and
the boundedness of $(U^0_{\2,h})_h$ in $L^1(\Omega^{'}_\delta)$, then
the space translates of $\bar U_{\2,\D}$ on $\Omega^{'}$ are estimated
uniformly for all sequence $\size(\D_m)_m$ tend to zero.\\
In the same way, we prove the space translate for $\alpha = \1$.
\end{proof}

We state the following lemma on time translate of $\bar
U_{\alpha,\D}$.

\begin{lemma}\label{lem:translater-time}$\left(\text{Time translate of } \bar U_{\alpha,\D}\right)$. Under the assumptions $({H}\ref{hyp:H1})-({H}\ref{hyp:H6})$ . Let $\D$ be a finite volume
  discretization of $\Omega\times (0,T)$ in the sense of Definition
  \ref{def:disc_QT} and let $p_{\alpha,\D}$ be a solution of
  \eqref{eq:p0}-\eqref{def_disc:pc}. Then, there exists a positive
  constant $C>0$ depending on $\Omega$, $T$ such that the following
  inequality hold:
  \begin{equation}\label{eq:translater_temps}
    \int_{\Omega \times
      (0,T-\tau)}\abs{\bar U_{\alpha,\D}(t+\tau,x)-\bar U_{\alpha,\D}(t,x)}^2\,dx \,dt \le \tilde{\omega}(\tau),
\end{equation}
for all $\tau\in (0,T)$. Here $\tilde{\omega}:\R^+\to\R^+$ is a modulus
of continuity, i.e. $\lim_{\tau\to 0}\tilde{\omega}(\tau)=0$. 
\end{lemma}

We state without proof the following lemma on time translate of $\bar
U_{\alpha,\D}$. Following the lemma \ref{borismostafa},  the proof is a direct consequence of the estimations \eqref{est:p_\2lobale} and \eqref{est:discrete_beta}, then we omit it.

\section{Study of the limit}\label{sec:limite}
\begin{prop}\label{prop:conv}
  Let $\hhm$ be a sequence of
  finite volume discretizations of $\Omega\times (0,T)$  such that $\lim_{m\to +\infty} {\rm
    size}(\hh)=0$. Then there exists subsequences, still denoted
  $(s_{\alpha,\hh})_\m$, $(p_{\alpha,\hh})_\m$ verify the following
  convergences
  \begin{align}   
    &\Vert U_{\alpha,\hh}-\bar U_{\alpha,\hh}\Vert_{L^1(\Omega ')}\longrightarrow 0,\label{conv:U-U_D}&&\\
    &U_{\alpha,\hh}\longrightarrow U_\alpha  && \text{ strongly in }L^p(Q_T) \text{ and a.e. in } Q_T \text{ for all } p\geq1,\label{conv:U}\\
    & \nabla\eh\B(s_{\2,\hh}) {\longrightarrow} \nabla \B(s_\2)&&
    \text{ weakly in } (L^2(Q_T))^\dm,
    \label{conv:beta}\\
    & \nabla\eh p\ch {\longrightarrow} \nabla p &&
    \text{ weakly in } (L^2(Q_T))^\dm, \label{conv:nabla_p}\\
    &s_{\alpha,\hh}\longrightarrow s_\alpha && \text{ almost everywhere in } 
    Q_T,\label{conv:s_alpha}\\
    &p_{\alpha,\hh} {\longrightarrow} p_\alpha && \text{ almost everywhere in }      Q_T. \label{conv:p_alpha}
  \end{align} 
Furthermore,
\begin{align}
\label{conv:s_12}
& 0\le s_\alpha \le 1 \text{ a.e. in } Q_T,\\
& U_\alpha=\phi \rho_\alpha(p_\alpha)s_\alpha \text { a.e. in } Q_T.
\label{iden=U}
\end{align}
\end{prop}
\begin{proof}
  For the first convergence \eqref{conv:U-U_D} it is useful to
  introduce the following inequality, for all $a,b\in\R$,
 $$\int_0^1 |\theta
 a + (1-\theta)b|\,d\theta\geq \frac{1}{2}(|a|+|b|).$$ Applying this
 inequality to $a=U_{\alpha,\hh}^{n+1}-U_{\alpha,\hh}^{n}$,
 $b=U_{\alpha,\hh}^{n}-U_{\alpha,\hh}^{n-1}$, from the definition of
 $\bar U_{\alpha,\hh}$ we deduce
$$ \int_0^T\int_{\Omega '}|U_{\alpha,\hh}(t,x)-\bar U_{\alpha,\hh}(t,x)| \dd x \dd t \leq
2\,\int_0^{T\!+\h}\int_{\Omega '} |\bar U_{\alpha,\hh}(t\!+\!\h,x)-\bar
U_{\alpha,\hh}(t,x)|\dd x \dd t.
$$ 
Since $\h$ tends to zero as $\size(\D_m)\to 0$, estimate
\eqref{eq:translater_temps} in Lemma \ref{lem:translater-time} implies
that the right-hand side of the above inequality converges to zero as
$\size(\D_m)$ tends to zero, and this established (\ref{conv:U-U_D}).\\
By the Riesz-Frechet-Kolmogorov compactness criterion, the relative
compactness of $(\bar U_{\alpha,\hh})_\m$ in $L^1(Q_T)$ is a
consequence of the Lemmas \ref{lem:translater-espace} and
\ref{lem:translater-time}. Now, the convergence \eqref{conv:U} in
$L^1(Q_T)$ and a.e in $Q_T$ becomes a consequence of
\eqref{conv:U-U_D}.  Due to the fact that $U_{\alpha,\hh}$ is bounded,
we establish the convergence in $L^1(Q_T)$. This ensures the following
strong convergences
\begin{align*}
  & \rho_\alpha(p_{\alpha,\hh})s_{\alpha,\hh} \longrightarrow l_\alpha
  \quad \text{ in $L^1(Q_T)$ and a.e. in $Q_T$ }.
\end{align*} 

Denote by $u_\alpha= \rho_\alpha(p_\alpha)s_\alpha$. Define the map
$\mathbb{H} : \R^+ \times \R^+ \mapsto \R^+ \times [0,\B(1)]$ defined
by
\begin{equation}
  \mathbb{H}(u_\2,u_\1) = (p,\B(s_\2))
\label{def:H}
\end{equation}
where $u_\alpha$ are solutions of the system
\begin{align*}
  & u_\2(p,\B(s_\2)) = \rho_\2(p-\bar{p}(\B
  ^{-1}(\B(s_\2))))\B ^{-1}(\B(s_\2))\\ & u_\1(p,\B(s_\2)) =
  \rho_\1(p-\tilde{p}(\B ^{-1}(\B(s_\2))))(1-\B ^{-1}(\B(s_\2)).
\end{align*}
Note that $\mathbb{H}$ is well defined as a diffeomorphism, since
\begin{eqnarray*}
  \frac{\partial u_\2 }{\partial p}&=&\rho_1^\prime(p-\bar{p}(\B ^{-1}(\B(s_1))))\B ^{-1}(\B(s_1))\ge0\\
  \frac{\partial u_\2 }{\partial \B}&=&  \rho_1^\prime(p-\bar{p}(\B ^{-1}(\B(s_1))))[ -\bar{p}^\prime(\B ^{-1}(\B(s_1))) ({\B ^{-1}}^\prime (\B(s_1)))]\B ^{-1}(\B(s_1)) \\&+& \rho_1(p-\bar{p}(\B ^{-1}(\B(s_1)))) {\B ^{-1}}^\prime (\B(s_1))\ge0\\
  \frac{\partial u_\1 }{\partial p}&=&-\rho_2^\prime(p-\tilde{p}(\B ^{-1}(\B(s_1))))(1-\B ^{-1}(\B(s_1)))\ge0\\
  \frac{\partial u_\1 }{\partial \B}&=& \rho_2^\prime(p-\tilde{p}(\B ^{-1}(\B(s_1))))  [-\tilde{p}^\prime(\B ^{-1}(\B(s_1))) ({\B ^{-1}}^\prime (\B(s_1)))] [ 1-\B ^{-1}(\B(s_1))]\\ &-& \rho_2(p-\tilde{p}(\B ^{-1}(\B(s_1)))){\B ^{-1}}^\prime (\B(s_1))\le0,\\
\end{eqnarray*}
and if one of the saturations is zero the other one is one, this
conserves that the jacobian determinant of the map $\mathbb{H}^{-1}$
is strictly negative.\\

As the map ${\mathbb{H}}$ defined in \eqref{def:H} is continuous, we
deduce
\begin{align*}
  & p\ch \longrightarrow p \quad \text{ a.e. in } Q_T,\\
  & \B(s_{l,\hh}) \longrightarrow \B^{*}\quad \text{ a.e. in } Q_T.
\end{align*} 
Then, as $\B^{-1}$ is continuous, we deduce
$$
s_{l,\hh} \longrightarrow s_l = \B^{-1}(\B^{*})\quad \text{ a.e. in }
  Q_T,
$$
and the convergences \eqref{conv:s_alpha} hold.

Consequently and due to the relationship between the pressure of each
phase and the global pressure defined in \eqref{def:pression_globale},
then the convergences \eqref{conv:p_alpha} hold  
\begin{align*}
  & p_{\alpha,\hh} \longrightarrow p_\alpha \quad \text{ a.e. in } Q_T.
\end{align*} 


It follows from Proposition \ref{prop:estimation_pression} that, the
sequence $(\Grad\ch p\ch)_\m$ is bounded in $(L^2(Q_T))^\dm$, and as a
consequence of the discrete Poincar\'e inequality, the sequence
$(p\ch)_\m$ is bounded in $L^2(Q_T)$.  Therefore there exist two
functions $p\in L^2(Q_T)$ and $\psi \in (L^2(Q_T))^\dm $ such that
\eqref{conv:nabla_p} holds and
$$\Grad\ch p\ch\longrightarrow \psi \text{ weakly in } (L^2(Q_T))^\dm.$$
It remains to identify $\nabla p$ by $\psi$ in the sense of
distributions. For that, it is enough to show as $m\to +\infty$:
$$ 
E_m:=\int \int_{Q_T}\nabla\ch p\ch \cdot\varphi \,\dd x \dd t
+\int \int_{Q_T} p\ch \di\varphi \,\dd x \dd t \longrightarrow
0,\quad\forall \varphi\in D(Q_T).
$$ Let $\hh$ be small enough such that
$\varphi$ vanishes in $T^{\text{ext}}_{K,\sigma}$ for all
$K\in\mathcal{T}$, then 
\begin{multline*}
  \int_{\Omega} p\ch \di \varphi(t,x) \,dx=\sum_{K\in \mathcal{T}}
  \int_K p\ch \di \varphi(t,x)\,dx\\
  =\sum_{K\in \mathcal{T}}\sum_{L\in N(K)} p^n_K
  \int_{\sigma_{\KL}}\varphi (t,x)\cdot \eta_{\KL}\,\dd \bord
  =\frac{1}{2}\sum_{K\in \mathcal{T}}\sum_{L\in
    N(K)}(p^n_K-p^n_L)\int_{\sigma_{\KL}}\varphi (t,x)\cdot
  \eta_{\KL}\,\dd \bord.
\end{multline*}
Now, from the definition of the discrete gradient,
\begin{align*}
  \int_{\Omega} \nabla\ch p\ch \varphi(t,x) \,dx
  &=\frac{1}{2}\sum_{K\in \mathcal{T}}\sum_{L\in N(K)}
  \int_{T_{\KL}}\nabla\ch p\ch  \varphi(t,x) \,dx\\
  &=\frac{1}{2}\sum_{K\in \mathcal{T}}\sum_{L\in
    N(K)}\frac{\dm}{\dis }(p^n_\Lll-p^n_K) \int_{T_{\KL}}
  \varphi(t,x)\cdot \eta_{\KL} \,dx
\end{align*}
Then,
\begin{equation*}
\begin{split}
  E_m= \frac{1}{2}\sum_{K\in \mathcal{T}}\sum_{L\in
    N(K)}\sigma_{\KL}(p^n_\Lll-p^n_\K) \Big(\frac{1}{|\sigma_{\KL}|}
  \int_{\sigma_{\KL}}\varphi (t,x)\cdot \eta_{\KL}\dd \bord  -
  \frac{1}{|T_{\KL}|}\int_{T_{\KL}} \varphi(t,x)\cdot \eta_{\KL} \,dx
  \Big)
\end{split}
\end{equation*}
Due to the smoothness of $\varphi$, one gets
$$ \Big| \frac{1}{|\sigma_{\KL}|} \int_{\sigma_{\KL}}
\varphi (t,x)\cdot \eta_{\KL}\dd \bord - \frac{1}{\abs{T_{\KL}}}
\int_{T_{\KL}} \varphi(t,x)\cdot \eta_{\KL} \,dx\Big| \leq C\; h,
$$
and the Cauchy-Scharwz inequality with the estimate \eqref{est:p_alpha} in Proposition
\ref{prop:estimation_pression} yield
\begin{align*}
  |E_m|\leq C h \sum_{n=0}^{N-1}\delta t \sum_{K\in
    \mathcal{T}}\sum_{L\in N(K)} |\sigma_{\KL}||p^n_L-p^n_K|
  \leq C h \sum_{n=0}^{N-1}\delta t \sum_{K\in
    \mathcal{T}}\sum_{L\in N(K)} |\sigma_{\KL}|\dis 
  \leq C h |\Omega| T.
\end{align*}
The identification of the limit in \eqref{iden=U} follows from the
previous convergence.
%
\end{proof}

\subsection{Proof of theorem \ref{theo:principal}}


  Let $T$ be a fixed positive constant and $\varphi\in D([0,T)\times
  \overline{\Omega})$. Set $\varphi_K^n:=\varphi(t^n,x_K)$ for all
  $K\in \mathcal{T}$ and $n\in[0,N]$.\\
  For the discrete liquid equation, we multiply the equation
  \eqref{sys_disc:pl} by $\pas \varphi_K^{n+1}$ and sum over $K\in
  \mathcal{T}$ and $n\in \{0,...,N\}$. This yields
  $$
  S_1^m +S_2^m+S_3^m + S_4^m = 0,
  $$
  where 
  \begin{equation*}
    \begin{split}
      S_1^m & = \snk \abs{K}\phi_K \left(
        \rho_\2(p^{n+1}_{\2,K})s^{n+1}_{\2,K}
        -\rho_\2(p^{n}_{\2,K})s^{n}_{\2,K}\right) \varphi_K^{n+1},\\
      S_2^m & = \sm \tokl \rho^{n+1}_{\2,\KL}G_\2(s^{n+1}_{\2,K},s^{n+1}_{\2,L};\de(p_{\2})) \varphi_K^{n+1},\\
      S_3^m & = \sm |\sigma_{\KL}|\Big((\rho^{n+1}_{\2,\KL})^2
      M_\2(s^{n+1}_{\2,K})({\bf g}_{\KL})^+ -
      (\rho^{n+1}_{\2,\KL})^2M_\2(s^{n+1}_{\2,L})({\bf
        g}_{\KL})^-\Big)\varphi_K^{n+1},\\
      S_4^m & = \sm \abs{K} \left( \rho_\2(p_{\2,K}^{n+1})
        s_{\2,K}^{n+1} f_{P,K}^{n+1}\varphi_K^{n+1}-
        \rho_\2(p_{\2,K}^{n+1}) \sli
        f_{I,K}^{n+1}\varphi_K^{n+1}\right).
    \end{split}
  \end{equation*}
  Making summation by parts in time and keeping in mind that
  $\varphi(T,x_K)=\varphi_K^{N+1}=0$. For all $K\in\mathcal{T}$, we
  get
  \begin{equation*}
    \begin{split}
      S_1^m = & - \snk
      \abs{K}\phi_K\rho_\2(p^{n+1}_{\2,K})s^{n+1}_{\2,K}\left(
        \varphi_K^{n+1}- \varphi_K^{n} \right) -\sum_{K\in
        \mathcal{T}_h}\abs{K}\phi_K\rho_\2(p^{0}_{\2,K})s^{0}_{\2,K}\varphi_K^0 \\
      = & -\snk \int_{t^n}^{t^{n+1}}\int_{K} \phi_K
      \rho_\2(p^{n+1}_{\2,K})s^{n+1}_{\2,K} \partial_t \varphi(t,x_K)\dd x\dd t 
       -\sum_{K\in
        \mathcal{T}_h}\int_{K}\phi_K\rho_\2(p^{0}_{\2,K})s^{0}_{\2,K}\varphi(0,x_K)
      \dd x.
    \end{split}
  \end{equation*}
Since $\phi\ch\rho_\2(p_{\2,\hh})s_{\2,\hh}$ and $\phi\ch
\rho_\2(p^0_{\2,\hh})s^0_{\2,\hh}$ converge almost everywhere
respectively to $\phi\rho_\2(p_\2)s_\2$ and $\phi
\rho_\2(p^0_\2)s^0_\2$, and as a consequence of Lebesgue dominated
convergence theorem, we get
\begin{equation*}
  \lim_{m\to +\infty}S_1^m =
  \int_{Q_T}\phi \rho_\2(p_\2)s_\2 \partial_t\varphi(t,x) \dd x \dd t 
  - \int_{\Omega}\phi \rho_\2(p^0_\2)s^0_\2 \varphi(0,x)\dd x.
\end{equation*}
Now, let us focus on convergence of the degenerate diffusive term to show 
\begin{equation}\label{conv:S2h}
  \begin{split}
    \lim_{m\to +\infty} S_2^m = -
    \int_{Q_T}\rho_\2(p_\2)M_\2(s_\2)\nabla p_\2\cdot \nabla \varphi
    \dd x \dd t.
  \end{split}
\end{equation}
Since the discrete gradient of each phase is not bounded, it is not possible to justify the pass to the limit in a straightforward way. To do this, we use the feature of global pressure and the auxiliary pressures defined in \eqref{def:pression_globale} and the discrete energy estimates in proposition \ref{prop:estimation_pression} and corollary \ref{cor:est}.

Gathering by edges, the term $S_2^m$ can be rewritten as:
  \begin{align*}
    S_2^m = & - \frac 1 2 \sm \tokl
    \rho^{n+1}_{\2,\KL}G_\2(s^{n+1}_{\2,K},s^{n+1}_{\2,L};\de(p_{\2}))
    \left( \varphi(t^{n+1},x_L) - \varphi(t^{n+1},x_K)\right)\\
    = \; & \frac 1 2 \sm \tokl
    \rho^{n+1}_{\2,\KL}M_\2(s^{n+1}_{\2,\KL}) \de(p_{\2}) \left(
      \varphi(t^{n+1},x_L) -
      \varphi(t^{n+1},x_K)\right) \\
    = \; & A_1^m+A_2^m,
  \end{align*}
with, by using the definition \eqref{def:pression_globale}, 
  \begin{align*}
    & A_1^m = \frac 1 2 \sm
    \tokl\rho^{n+1}_{\2,\KL}M_\2(s^{n+1}_{\2,\KL}) \de(p)\left(
      \varphi(t^{n+1},x_L) -
      \varphi(t^{n+1},x_K)\right), \\
    & A_2^m = - \frac 1 2 \sm \tokl
    \rho^{n+1}_{\2,\KL}M_\2(s^{n+1}_{\2,\KL}) \de(\bar{p}(s_\3))\left(
      \varphi(t^{n+1},x_L) - \varphi(t^{n+1},x_K)\right).
  \end{align*}
Let us show that 
\begin{align}\label{conv:A1}
  \lim_{m\to +\infty} A_1^m = \int_{Q_T} \rho_\2(p_\2)M_\2(s_\2)
  \nabla p \cdot \nabla \varphi\, \dd x \dd t.
\end{align}
For each couple of neighbours $K$ and $L$ we denote $\slmin$ the minimum of
$s_{\2,K}^{n+1}$ and $s_{\2,L}^{n+1}$ and we introduce
\begin{align*}
  A_1^{m,*} = \frac 1 2 \sm \tokl\rho^{n+1}_{\2,\KL}M_\2(\slmin)
  \de(p)\left( \varphi(t^{n+1},x_L) - \varphi(t^{n+1},x_K)\right)
\end{align*}
Remark that
\begin{align*}
  A_1^{m,*} & = \frac 1 2 \sm \dm |T_{\KL}| 
  \rho^{n+1}_{\2,\KL}M_\2(\slmin) \frac{p_L - p_K
  }{\dis }\frac{\varphi(t^{n+1},x_L) -
    \varphi(t^{n+1},x_K)}{\dis }\\ & = \frac 1 2 \sm
  |T_{\KL}|\rho^{n+1}_{\2,\KL}M_\2(\slmin) \nabla_{\KL} p_{\D_m}\cdot
  \eta_{\KL}\nabla \varphi(t^{n+1},x_{\KL})\cdot \eta_{\KL},
\end{align*}
where $x_{\KL} = \theta x_K + (1-\theta)x_L$, $0<\theta <1$, is some
point on the segment $]x_K,x_L[$. Recall that the value of
$\nabla_{\KL}$ is directed by $\eta_{\KL}$, so
$$
\nabla_{\KL} p_{\D_m}\cdot \eta_{\KL}\nabla
\varphi(t^{n+1},x_{\KL})\cdot \eta_{\KL}= \nabla_{\KL} p_{\D_m}\cdot
\nabla \varphi(t^{n+1},x_{\KL})
$$
Define $\overline{s}_{\alpha,\D_m}$ and $\underline{s}_{\alpha,\D_m}$
by
$$
\overline{s}_{\alpha,\D_m}|_{(t^n,t^{n+1}]\times
  T_{\KL}}:=\max\{s_{\alpha,K},s_{\alpha,L}\},\quad
\underline{s}_{\alpha,\D_m}|_{(t^n,t^{n+1}]\times
  T_{\KL}}:=\min\{s_{\alpha,K},s_{\alpha,L}\}
$$
Now, $A_1^{m,*}$ can be written under the following continues form
$$A_1^{m,*} = \int_{0}^{T}\int_{\Omega} \rho_\2(p_{\2,\D_m}) 
M_\2(\underline{s}_{\2,\D_m}) \nabla_{\D_m}p_{\D_m}\cdot (\nabla
\varphi)_{\D_m}\dd x \dd t.$$ By the monotonicity of $\B$ and thanks
to the estimate \eqref{est:discrete_beta}, we have
\begin{align*}
  \int_{0}^{T}\int_{\Omega}
  \abs{\B(\overline{s}_{\2,\D_m})-\B(\underline{s}_{\2,\D_m})}^2\dd x \dd
  t \le & \sm |T_{\KL}| \left(\B(s_{\3,L}^{n+1}) - \B(s_{\3,K}^{n+1})
  \right)^2\\ &\le C \size(\T)^2\sm
  \frac{\abs{\sig}}{\dis}\abs{\B(s_{\3,L}^{n+1}) -
    \B(s_{\3,K}^{n+1})}^2\\ & \le C \size(\T)^2.
\end{align*}
Since $\B^{-1}$ is continuous, we deduce up to a subsequence
\begin{align}\label{conv:sminsmax}
  \abs{\underline{s}_{\alpha,\D_m}-\overline{s}_{\alpha,\D_m}}\to 0
  \text{ a.e. on } Q_T.
\end{align}
Moreover, we have $\underline{s}_{\alpha,\D_m}\le
s_{\alpha,\D_m}\le\overline{s}_{\alpha,\D_m} $ and $s_{\alpha,\D_m}\to
s_\alpha$ a.e. on $Q_T$. Consequently, and due to the continuity of
the mobility function $M_\2$ we have $M_\2(\underline{s}_{\2,\D_m})\to
M_\2(s_\2)$ a.e on $Q_T$ and in $L^p(Q_T)$ for $p<+\infty$.

As consequence of the convergence \eqref{conv:p_alpha} and by the
Lebesgue dominated convergence theorem we get
$$
\rho_\2(p_{\2,\D_m})M_\2(\underline{s}_{\2,\D_m}) (\nabla
\varphi)_{\D_m} \to \rho_\2(p_\2) M_\2(s_\2)\nabla \varphi \text{
  strongly in } (L^{2}(Q_T))^\dm.
$$
And as consequence of the weak convergence on global pressure \eqref{conv:nabla_p}, we obtain that
$$
  \lim_{m\to +\infty} A_1^{m,*} = \int_{Q_T} \rho_\2(p_\2)M_\2(s_\2)
  \nabla p \cdot \nabla \varphi \, \dd x \dd t.
$$
It remains to show that 
\begin{equation}\label{A1m-A1m*}
  \lim_{m\to +\infty} \abs{A_1^m-A_1^{m,*}}=0.
\end{equation}
Remark that 
\begin{equation*}
  \abs{M_\2(s_{\2,\KL}^{n+1})\de(p)-M_\2(\slmin)\de(p)}\le
  C \abs{s_{\2,L}^{n+1}-s_{\2,K}^{n+1}}\abs{\de(p)}.
\end{equation*}
Consequently
$$
  \abs{A_1^m-A_1^{m,*}} \le C\int_{Q_T} 
  \abs{s_{\2,L}^{n+1}-s_{\2,K}^{n+1}} \nabla_{\D_m}p_{\D_m}\cdot
  (\nabla \varphi)_{\D_m}\dd x \dd t.
$$
Applying the Cauchy-Schwarz inequality, and thanks to the uniform
bound on $\nabla_{\D_m}p_{\D_m}$ and the convergence
\eqref{conv:sminsmax}, we establish \eqref{A1m-A1m*}.\\

To prove the pass to limit of $A_2^m$, we need to prove firstly that
\begin{align*}
  \| \de(\Gamma(s_\3)) - \sqrt{M_\2(s^{n+1}_{\2,\KL})}
  \de(\bar{p}(s_\3))\|_{L^2(Q_T)} \to 0 \text{ as }
  \text{size}(\T) \to 0,
\end{align*}
where $\Gamma(s_\3)=\int_{0}^{s_\3}\sqrt{M_\2(z)}\frac{\dd
  \bar{p}}{\dd s_\3}(z) \dd z$.\\
In fact.  Remark that there exist $a\in[s_{\3,K},s_{\3,L}]$ such as:
 \begin{align*}
   |\de(\Gamma(s_\3)) - \sqrt{M_\2(s^{n+1}_{\2,K,L})}\de(\bar{p}(s_\3))|& 
   = |\sqrt{M_\2(a)} - \sqrt{M_\2(s^{n+1}_{\2,K,L})}|| 
        \de(\bar{p}(s_\3))|\\ & 
   \le C |\de(\bar{p}(s_\3))| 
   \le C \abs{s_{\3,L}^{n+1}-s_{\3,K}^{n+1}} \\ &
   \le C \abs{\B(s_{\3,L}^{n+1}) - \B(s_{\3,K}^{n+1})}^\theta,
    \end{align*}
since $\B^{-1}$ is an H\"older function. Thus we get,
 \begin{equation*}
   \begin{aligned}
     & \| \de(\Gamma(s_\2)) - \sqrt{M_\2(s^{n+1}_{\2,K,L})}
     \de(\bar{p}(s_\3))\|^2_{L^2(Q_T)}\\ &
          =  \sm |T_{K,L}| | \de(\Gamma(s_\2)) - \sqrt{M_\2(s^{n+1}_{\2,K,L})}
           \de(\bar{p}(s_\3))|^2 \\ & 
     \le \sm |T_{K,L}|^{1-\theta} |T_{K,L}|^{\theta}
         \abs{\B(s_{\3,L}^{n+1}) - \B(s_{\3,K}^{n+1})}^{2\theta},
       \end{aligned}
 \end{equation*}
 and using the Cauchy-Schwarz inequality and the estimate , we deduce 
  \begin{equation*}
   \begin{aligned}
     & \| \de(\Gamma(s_\2)) - \sqrt{M_\2(s^{n+1}_{\2,K,L})}
     \de(\bar{p}(s_\3))\|^2_{L^2(Q_T)}\\ &
       \le \left(\sm |T_{K,L}|\right)^{1-\theta}\left(\sm |T_{K,L}|
         \abs{\B(s_{\3,L}^{n+1}) - \B(s_{\3,K}^{n+1})}^2\right)^{\theta}\\ &
    \le C (\size(\T))^{2\theta} \left(\sm \frac{\abs{\sigma_{K,L}}}{\dis}
            \abs{\B(s_{\3,L}^{n+1}) - \B(s_{\3,K}^{n+1})}\right)^\theta 
       \end{aligned}
 \end{equation*}
 which shows that $ \| \de(\Gamma(s_\2)) -
 \sqrt{M_\2(s^{n+1}_{\2,\KL})} \de(\bar{p}(s_\3))\|^2_{L^2(Q_T)} \to
 0$ as $\text{size}(\T) \to 0$.  And from \eqref{est:discrete_pbar} in corollary \ref{cor:est},
 we deduce that there exists a constant $C>0$ where the following
 inequalities hold:
  \begin{align}\label{est:discrete_Gamma}
    \sm \tokl(\de(\Gamma(s_\3)))^2\le C.
  \end{align}
  That prove 
\begin{align}\label{est:Gamma}
  \nabla_{\D_m} \Gamma(s_{\3,\D_m})\to \nabla \Gamma (s_\3) \text{
    weakly in } (L^2(Q_T))^\dm.
\end{align}
 As consequence
\begin{align}\label{est:sqrtbarp}
  \sqrt{M_\2(s_{\2,\D_m})}\nabla_{\D_m} \bar{p}(s_{\3,\D_m})\to
  \nabla \Gamma (s_\3) \text{ weakly in } (L^2(Q_T))^\dm.
\end{align}
Rearranging $A_2^{m}$ to write 
\begin{equation*}
    A_2^{m}  = - \frac 1 2 \sm
    \abs{T_{\KL}} \rho^{n+1}_{\2,\KL}M_\2(s^{n+1}_{\2,\KL})
    \nabla_{\KL} \bar{p}(s_{\3,\D_m})\cdot \eta_{\KL}\nabla
    \varphi(t^{n+1},x_{\KL})\cdot \eta_{\KL},
\end{equation*}
where $x_{\KL} = \theta x_K + (1-\theta)x_L$, $0<\theta <1$, is some
point on the segment $]x_K,x_L[$. 
using again that the mesh is orthogonal, we can write 
$$A_2^{m} = - \int_{0}^{T}\int_{\Omega} \rho_\2(p_{\2,\D_m}) 
M_\2(s_{\2,\D_m}) \nabla_{\D_m}\bar{p}(s_{\3,\D_m})\cdot (\nabla
\varphi)_{\D_m}\dd x \dd t.
$$ 
As a consequence of the convergences \eqref{conv:s_alpha},
\eqref{conv:p_alpha} and by the Lebesgue theorem we get
$$
\rho_\2(p_{\2,\D_m}) \sqrt{M_\2(s_{\2,\D_m})} (\nabla \varphi)_{\D_m}
\to \rho_\2(p_\2) \sqrt{M_\2(s_\2)}\nabla \varphi \text{ strongly in }
(L^{2}(Q_T))^\dm.
$$
And as consequence of \eqref{est:sqrtbarp},
\begin{align}\label{lim:A2m}
  \lim_{m\to +\infty} A_2^{m} & = - \int_{0}^{T}\int_{\Omega}
  \rho_\2(p_{\2}) \sqrt{M_\2(s_{\2})} \nabla \Gamma(s_{\3})\cdot
  \nabla \varphi\dd x \dd t\\ & = - \int_{0}^{T}\int_{\Omega}
  \rho_\2(p_{\2}) M_\2(s_{\2}) \nabla \bar{p}(s_{\3})\cdot \nabla
  \varphi\dd x \dd t.
\end{align}
Now, we treat the convergence of the gravity term
\begin{align}\label{conv:gravityg}
  \lim_{m\to +\infty} S_3^m = - \int_{0}^{T}
  \int_{\Omega}\rho_\2(p_\2)M_\2(s_\2)\G\cdot \nabla \varphi \dd x \dd
  t.
\end{align}
Perform integration by parts \eqref{integration_part}
\begin{align*}
  S_3^m & = \sm F_{\2,\KL}^{n+1} \varphi(t^{n+1},x_K)\\ & = -\frac 1 2
  \sm F_{\2,\KL}^{n+1} \left( \varphi(t^{n+1},x_L) -
    \varphi(t^{n+1},x_K) \right).
\end{align*}
Note that the numerical flux  $F_{\2,\KL}^{n+1} $ is independent of the gradient of pressures and the pass to the limit on $S_3^m$ is mush simple then the term $A_1^{m,*}$ since the discrete gradient of global pressure is replaced by the gravity vector {\bf g}. We omit this proof of \eqref{conv:gravityg}.

Finally, $S_4^m$ can be written equivalently
\begin{align*}
  S_4^m = & \snk\int_{t^n}^{t^{n+1}} \int_{K}
  \rho_\2(p_{\2,K}^{n+1})s_{\2,K}^{n+1}f_P(t,x)
  \varphi(t^{n+1},x_K)\dd x \dd t \\ & \hspace{-.5cm}- \snk\int_{t^n}^{t^{n+1}}
  \int_{K} \rho_\2(p_{\2,K}^{n+1}) \sli
  f_{I}(t,x)\varphi(t^{n+1},x_K)\dd x \dd t.
\end{align*}
From the convergences \eqref{conv:s_alpha}, \eqref{conv:p_alpha} and
by the Lebesgue dominated convergence theorem, we get
\begin{align*}
\lim_{m\to +\infty} S_4^{m,*} =  \int_{Q_T} \rho_\2(p_\2)s_\2 f_P(t,x)\varphi(t,x)\dd x \dd t
 \hspace{-.5cm} - \int_{Q_T} \rho_\2(p_\2)s_\2^I f_I(t,x)\varphi(t,x)\dd x \dd t,
\end{align*}
which completes the proof of the theorem \ref{theo:principal}. 



 \section{Numerical results}\label{sec:numerical}
 

 In this section we show some numerical experiments simulating the five spot problem in petroleum engineering.  A Newton algorithm is implemented to approach the
 solution of nonlinear system \eqref{sys_disc:pl}-\eqref{sys_disc:pg}
 coupled with a bigradient method to solve linear system arising from
 the Newton algorithm process.

 We will provide two tests made on a nonuniform admissible grid.

 \begin{figure}[ht]
 \begin{center}
 \includegraphics[width=0.6\linewidth]{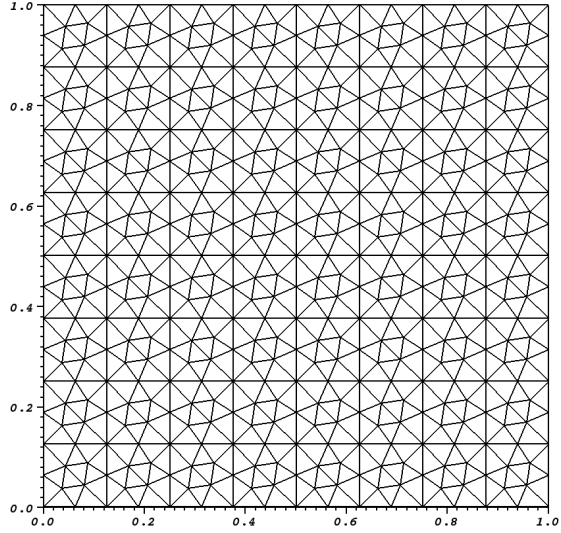} 
 \end {center}
 \caption{\footnotesize
 Mesh with $896$ triangles 
 }
 \label{fig:meshtriangle}
 \end{figure}

 \bigskip
 \vspace{3mm}
 Datas used for the numerical tests are the following :
 \[
 \begin{array}{ll}

 k_1(s_1) = s_{1}^{2},\, k_2(s_2)=s_2^2\\
 {\bf K}=0.15 10^{-10} \mbox{m}^{2},~\phi = 0.206,\\
 \mu_2=10^{-3} \mbox{ Pa.s}\text{(water viscosity)},~\mu_1=9 10^{-5}\mbox{ Pa.s}\text{(gas viscosity)},\\
 \rho(p) = \rho_{ref}(1+c_{ref}(p - p_{ref})),\mbox{ with }\rho_{ref} = 400 \text{ Kg}\,\text{m}^{-3},\,c_{ref}=10^{-6}\mbox{Pa}^{-1},\,p_{ref}=1.013\,10^{5} \text{ Pa},\\
 L_x =1 \text{m},\, L_y = 1\text {m}\text{ (the length and the width of the domain)}\\
 P_c(s)=P_{max}(1-s), \text{ with } P_{max}=10^{5}\text{Pa}.
 \end{array}
 \]

 {\bf Initial conditions. } Initially the saturation of gas is considered to be equal to $0.9$ in the whole domain and the gas pressure is considered to be $1.013\,10^5$ Pa.

 {\bf Boundary conditions.} The wetting fluid (water) is injected in the left-down corner  in the region $([0, 0.1]\times\{0\})\cup (\{0\}\times[0,0.1])$ 
 with a constant pressure equal to $4.026\,10^5$ Pa. 
 The right-top corner where  $([0.9, 1]\times\{1\})\cup (\{1\}\times[0.9,1])$
  keeps fluids flow freely at atmospheric pressure where as the rest of the boundary is 
 assumed to be impervious (zero fluxes  are imposed). The influence of boundary conditions can be seen in all figures.

 {\bf Meshes. } The domain is recovered   by $896$ admissible triangles see figure \ref{fig:meshtriangle}. 

\bigskip

\begin{figure}
 \begin{center}
  \begin{tabular}{cc}
 \includegraphics[width=0.4\linewidth]{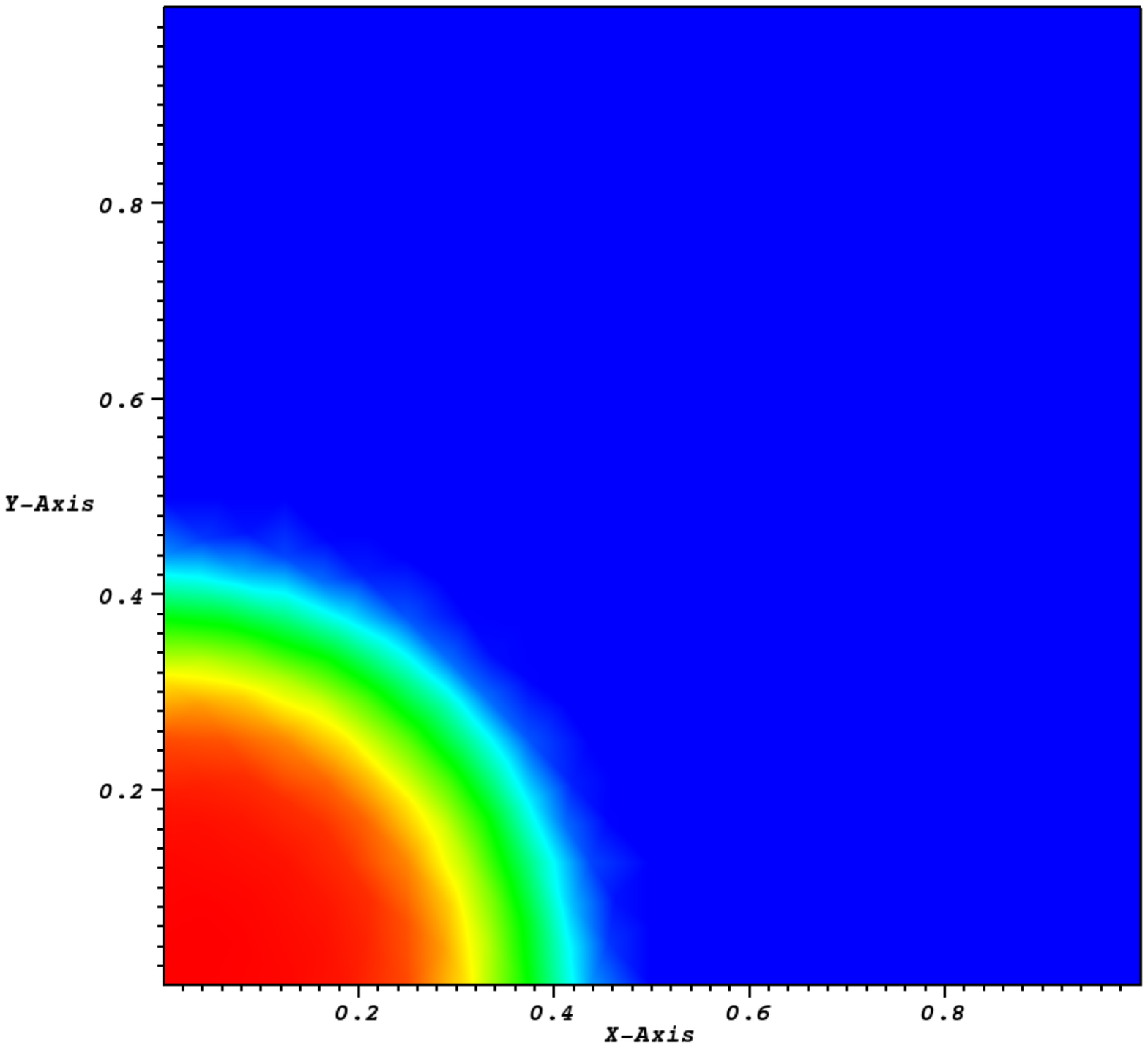} 
 &
 \includegraphics[width=0.4\linewidth]{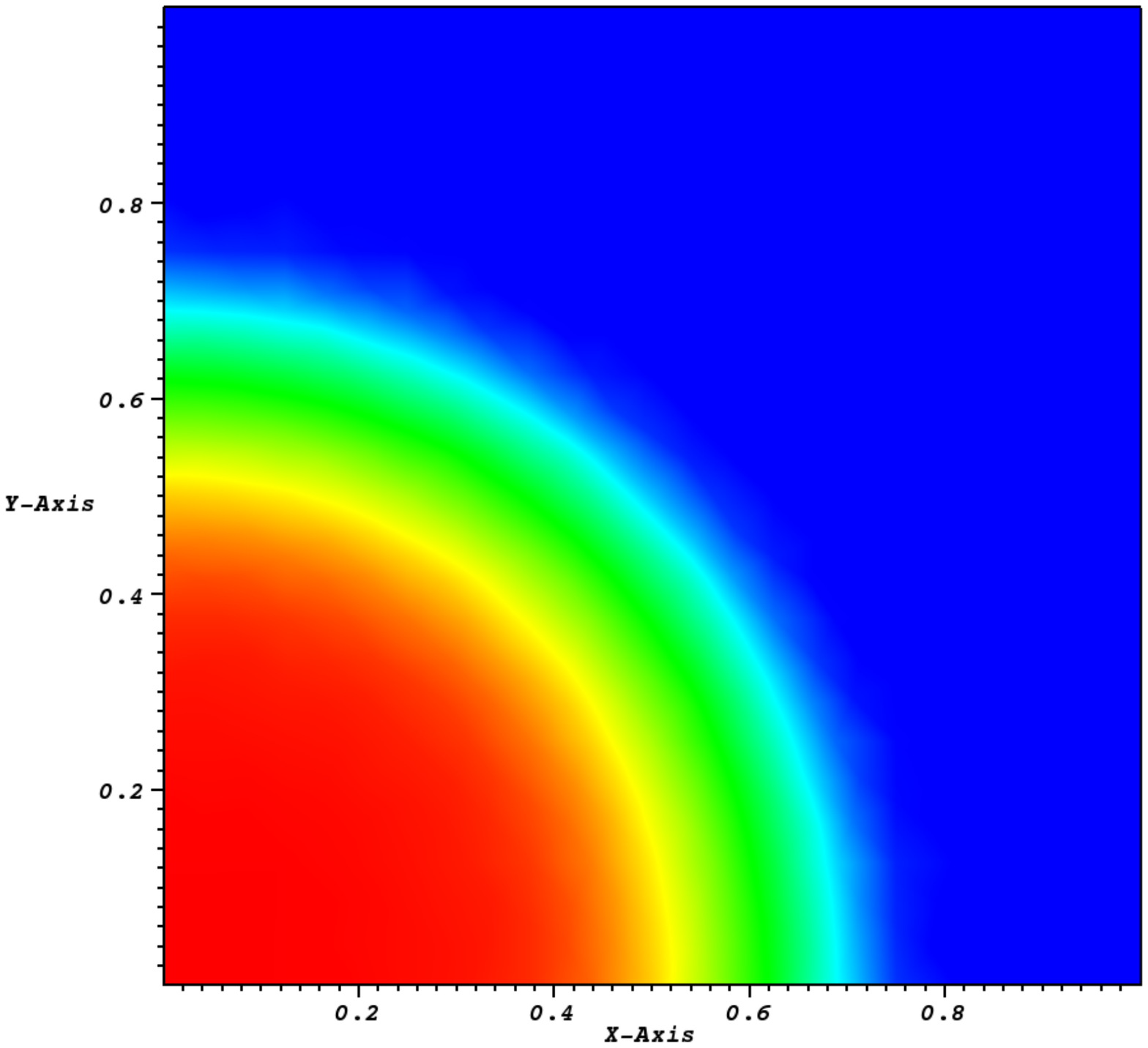} 
 \end{tabular}
 \caption{\footnotesize
 Water field including capillary effect at time $T=6$s (left) and  at time $T=20$s with $0.1\le s\le 1$.
 }
 \label{fig:test1:1}
 \end {center}
 \end{figure}

\begin{figure}
 \begin{center}
  \begin{tabular}{cc}
 \includegraphics[width=0.4\linewidth]{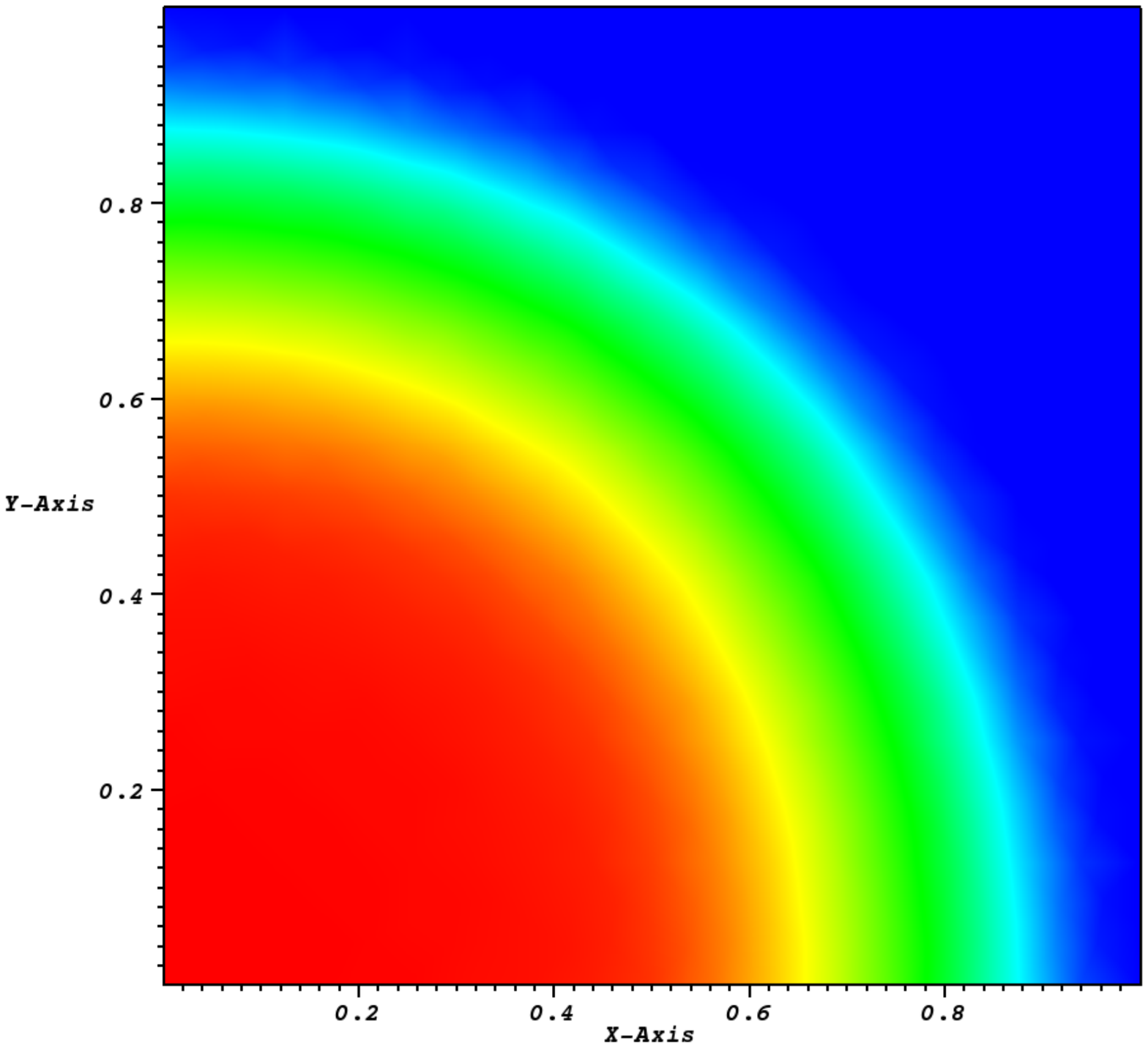} 
 &
 \includegraphics[width=0.4\linewidth]{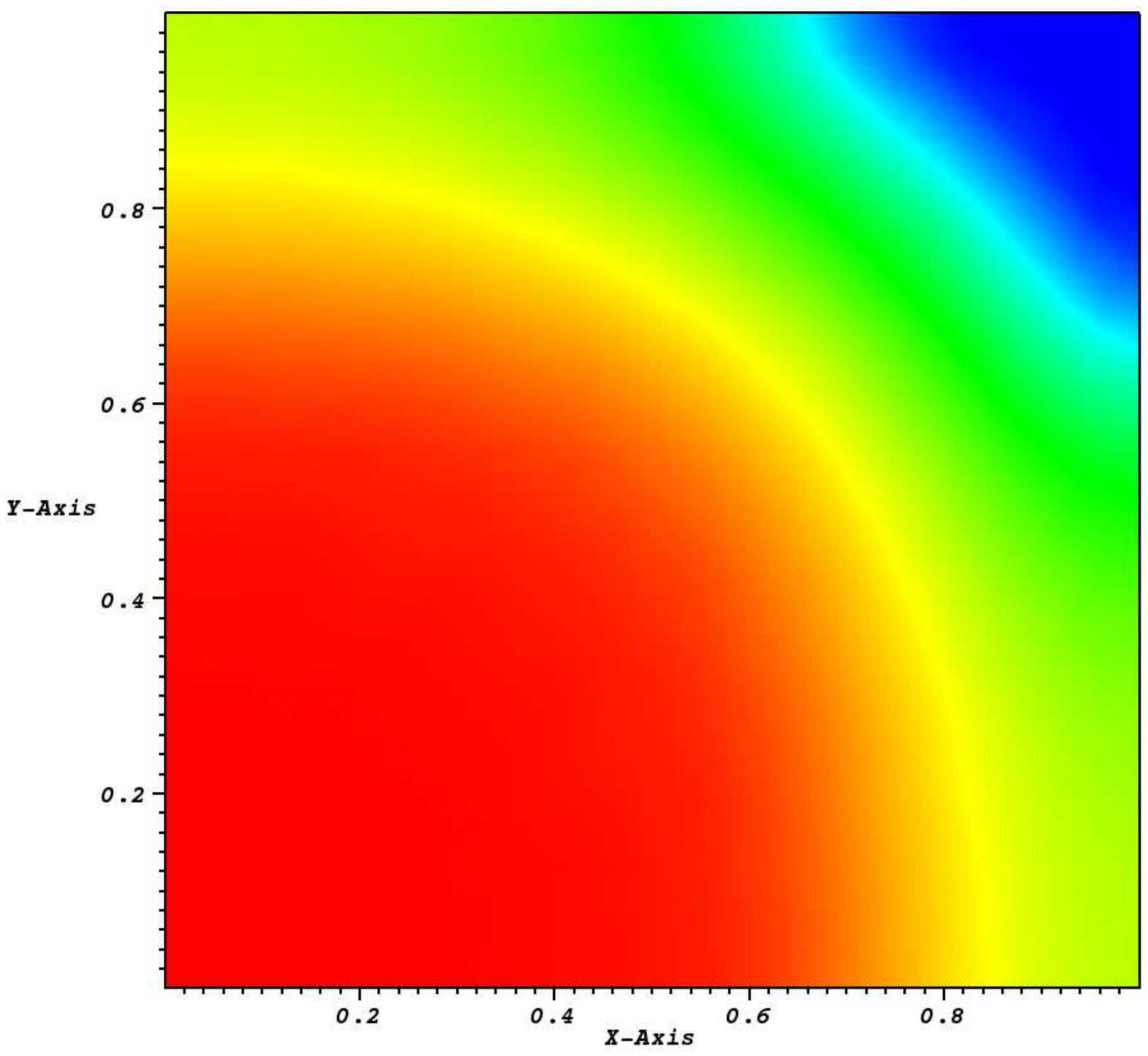} 
 \end{tabular}
 \caption{\footnotesize
 Water field including capillary effect at time $T=35$s (left) and  at time $T=60$s with $0.1\le s\le 1$.
 }
 \label{fig:test1:2}
 \end {center}
 \end{figure}

\begin{figure}
 \begin{center}
  \begin{tabular}{cc}
 \includegraphics[width=0.4\linewidth]{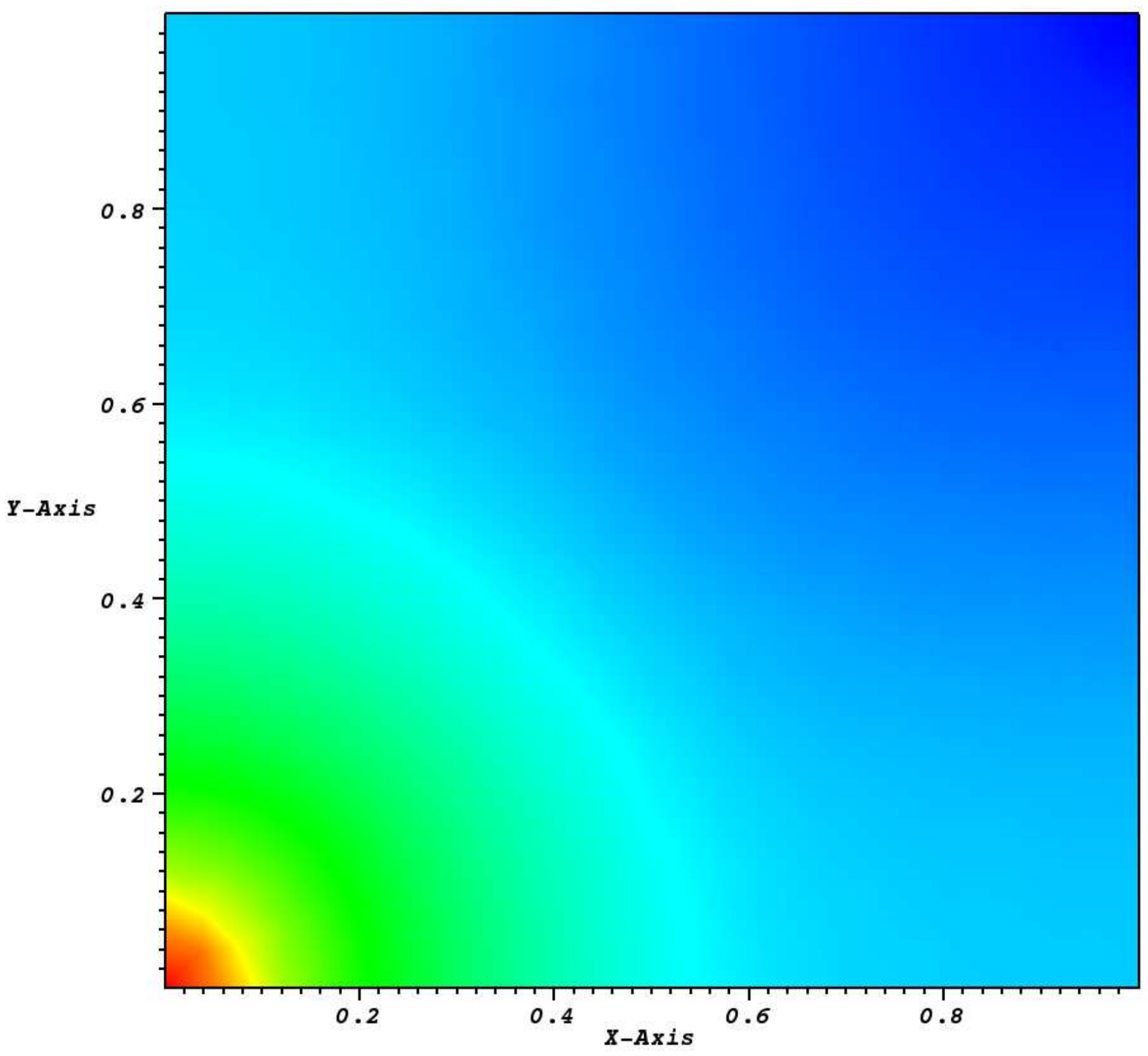} 
 &
 \includegraphics[width=0.4\linewidth]{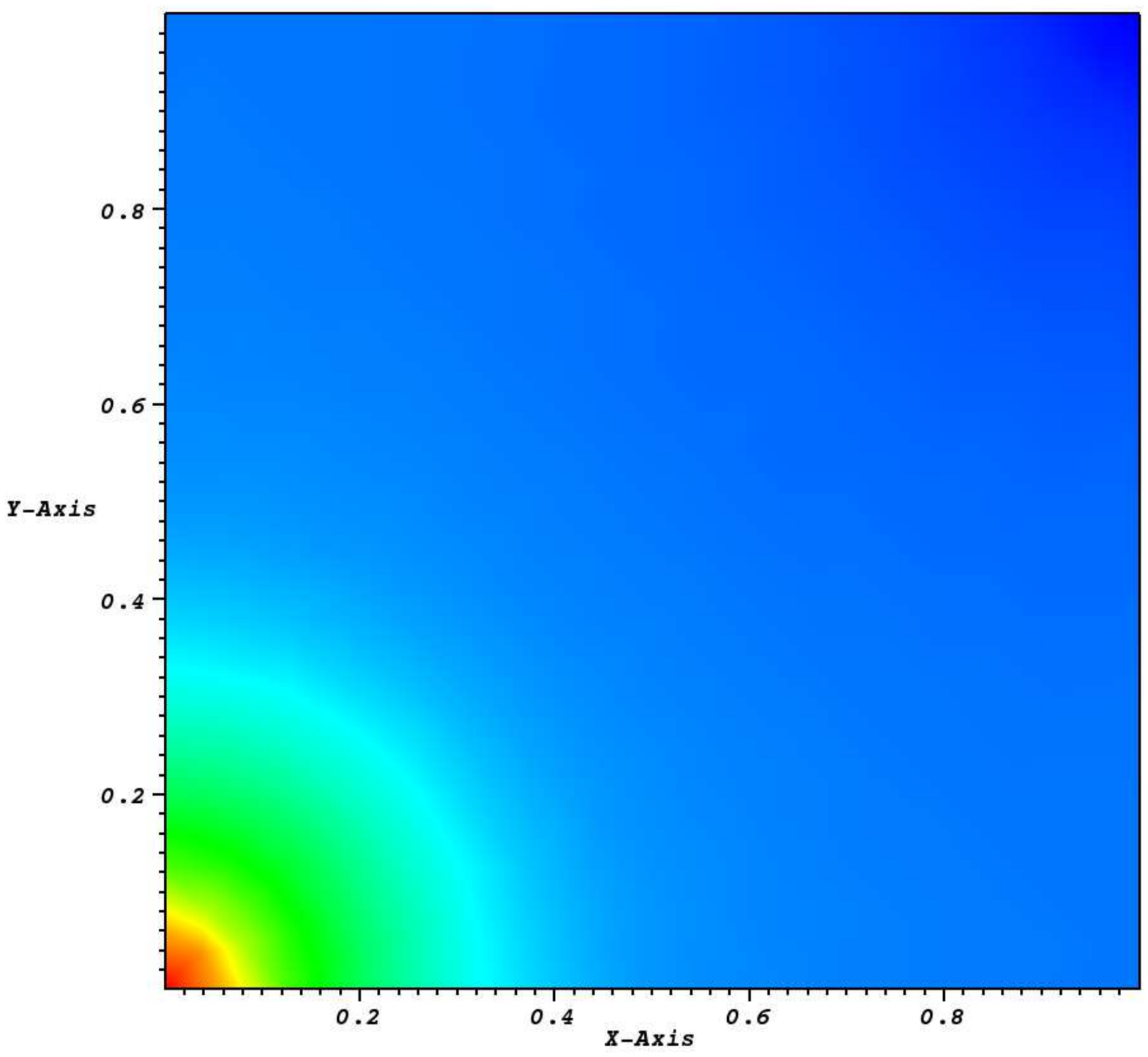} 
 \end{tabular}
 \caption{\footnotesize
 Gas pressure  field  including capillary effect at time $T=35$ (left) and  at time $T=60$ with $1.013\,\times 10^5\le P\le 4 \times 10^5$ Pa.
 }
 \label{fig:test1:3}
 \end {center}
 \end{figure}

\begin{figure}
 \begin{center}
  \begin{tabular}{cc}
 \includegraphics[width=0.4\linewidth]{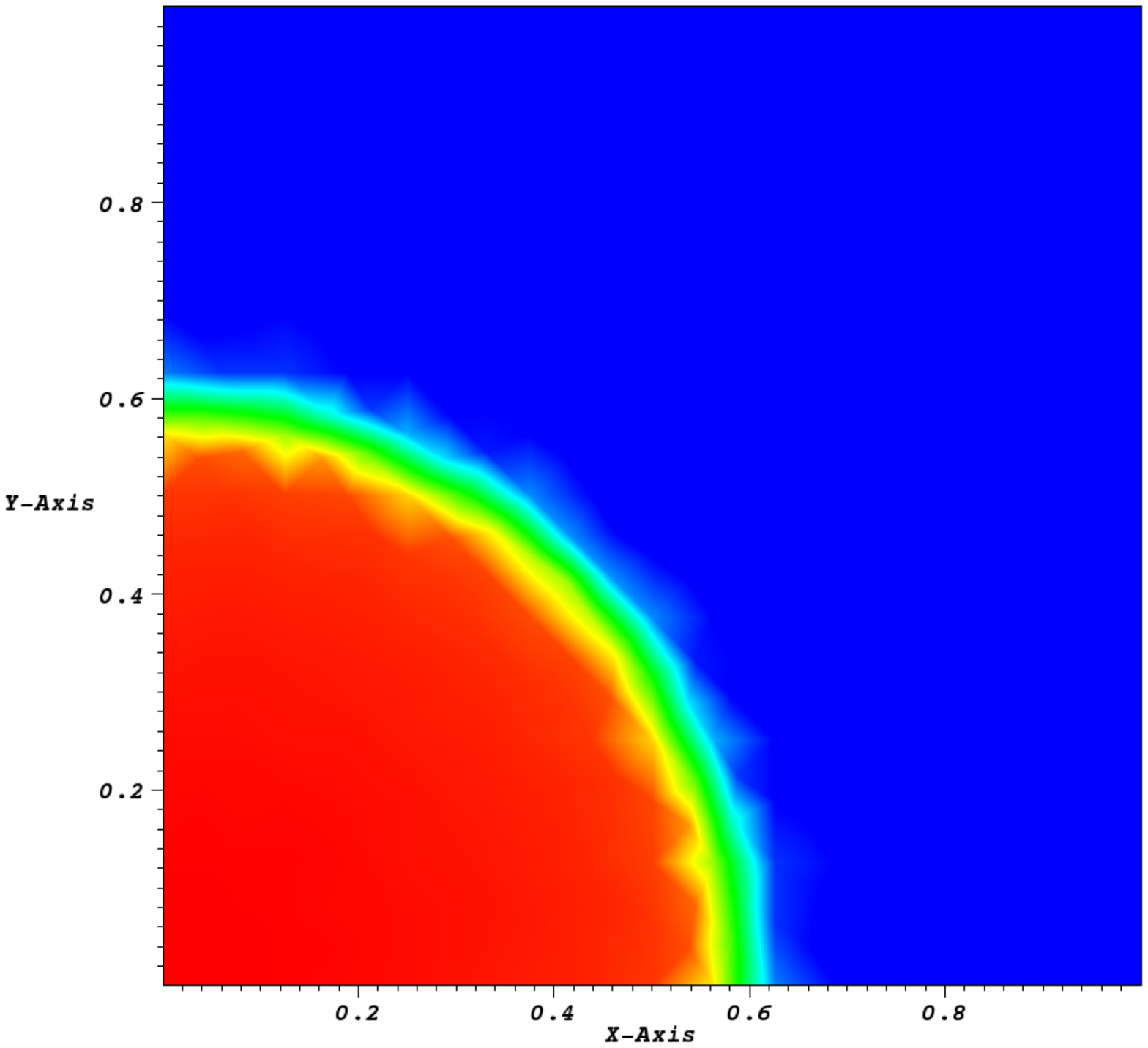} 
 &
 \includegraphics[width=0.4\linewidth]{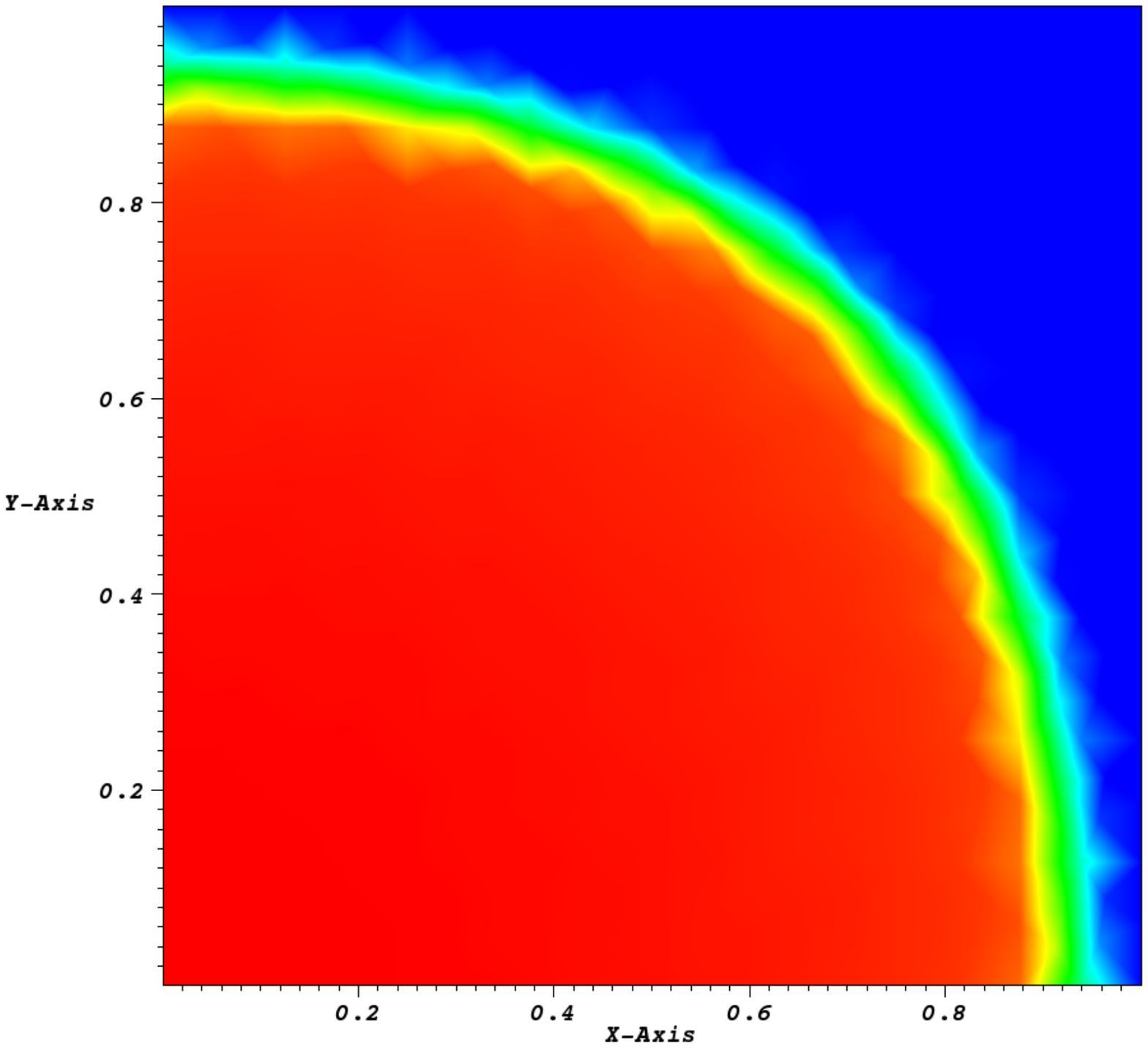} 
 \end{tabular}
 \caption{\footnotesize
 Water field without capillary terms at time $T=35$ (left) and  at time $T=60$ with $0.1\le s\le 1$
 }
 \label{fig:test1:4}
 \end {center}
 \end{figure}

 Figures \ref{fig:test1:1} - \ref{fig:test1:4} show the diffusive effects of the capillary terms, 
 notably the dissipation of chocs due to the hyperbolic operator Fig. \ref{fig:test1:4}. In fact, during the 
 stage of the displacement saturation shock propagate through rock for flows where 
 capillarity terms are neglected, see figure \ref{fig:test1:4}. This shock, where capillarity effects are 
 signifiant, it is diffused. However, a part of the the shock wave maintains its 
 sharp front.




\begin{thebibliography}{10}

\bibitem{Alt83}
H.W. Alt and S.~Luckhaus.
\newblock Quasilinear elliptic-parabolic differential equations.
\newblock {\em Math. Z., 3}, pages 311--341, 1983.

\bibitem{amaziane-ossmani}
B.~Amaziane and M.~El Ossmani.
\newblock Convergence analysis of an approximation to miscible fluid flows in
  porous media by combining mixed finite element and finite volume methods.
\newblock {\em Wiley InterScience (www.interscience.wiley.com). DOI
  10.1002/num.2029}, 2007.

\bibitem{amirat-fv}
Y.~Amirat, D.~Bates, and A.~Ziani.
\newblock Convergence of a mixed finite element-finite volume scheme for a
  parabolic-hyperbolic system modeling a compressible miscible flow in porous
  media.
\newblock {\em Numer. Math.}, 2005.

\bibitem{arbogast}
T.~Arbogast.
\newblock Two-phase incompressible flow in a porous medium with various non
  homogeneous boundary conditions.
\newblock {\em IMA Preprint series 606}, 1990.

\bibitem{aziz}
K.~Aziz and A.~Settari.
\newblock {\em Petroleum reservoir simulation}.
\newblock Applied Science Publishers LTD, London, 1979.

\bibitem{Bear67}
J.~Bear.
\newblock {\em Dynamic of flow in porous media}.
\newblock Dover, 1986.

\bibitem{ZS11fv}
M.~Bendahmane, Z.~Khalil, and M.~Saad.
\newblock Convergence of a finite volume scheme for gas water flow in a
  multi-dimensional porous media.
\newblock {\em submitted}, 2010.

\bibitem{brenier}
Y.~Brenier and J.~Jaffr\'e.
\newblock Upstream differencing for multiphase flow in reservoir simulation.
\newblock {\em SIAM J. Numer. Anal.}, 28:685--696, 1991.

\bibitem{CS10}
F.~Caro, B.~Saad, and M.~Saad.
\newblock Two-component two-compressible flow in a porous medium.
\newblock {\em Acta Applicandae Mathematicae (accepted)}, DOI:
  10.1007/s10440-011-9648-0 (2011).

\bibitem{chavent}
G.~Chavent and J.~Jaffré.
\newblock {\em Mathematical models and finite elements for reservoir
  simulation: single phase, multiphase, and multicomponent flows through porous
  media}.
\newblock North Holland, 1986.

\bibitem{chen99}
Z.~Chen.
\newblock Degenerate two-phase incompressible flow. existence, uniqueness and
  regularity of a weak solution.
\newblock {\em Journal of Differential Equations}, 171:203--232, 2001.

\bibitem{chen2002}
Z.~Chen.
\newblock Degenerate two-phase incompressible flow. regularity, stability and
  stabilization.
\newblock {\em Journal of Differential Equations}, 186:345--376, 2002.

\bibitem{chen}
Z.~Chen and R.~E. Ewing.
\newblock Mathematical analysis for reservoirs models.
\newblock {\em SIAM J. math. Anal.}, 30:431--452, 1999.

\bibitem{evans:book}
L.~Evans.
\newblock {\em Partial Differential Equations}.
\newblock American Mathematical Society, 2010.

\bibitem{Eymard:book}
R.~Eymard, T.~Gallou{\"e}t, and R.~Herbin.
\newblock {\em Finite Volume Methods}, volume~7.
\newblock Handbook of Numerical Analysis, P. Ciarlet, J. L. Lions, eds,
  North-Holland, Amsterdam, 2000.

\bibitem{Eymard00}
R.~Eymard, R.~Herbin, and A.~Michel.
\newblock Mathematical study of a petroleum-engineering scheme.
\newblock {\em Mathematical Modelling and Numerical Analysis}, 37(6):937--972,
  2003.

\bibitem{feng}
X.~Feng.
\newblock On existence and uniqueness results for a coupled systems modelling
  miscible displacement in porous media.
\newblock {\em J. Math. Anal. Appl.}, 194(3):883--910, 1995.

\bibitem{GM96}
G.~Gagneux and M.~Madaune-Tort.
\newblock {\em Analyse mathematique de models non lineaires de l'ingeniere
  petroli\`ere}, volume~22.
\newblock Springer-Verlag, 1996.

\bibitem{CS04}
C.~Galusinski and M.~Saad.
\newblock On a degenerate parabolic system for compressible, immiscible,
  two-phase flows in porous media.
\newblock {\em Advances in Diff. Eq.}, 9(11-12):1235--1278, 2004.

\bibitem{CS08}
C.~Galusinski and M.~Saad.
\newblock A nonlinear degenerate system modeling water-gas in reservoir flows.
\newblock {\em Discrete and Continuous Dynamical System}, 9(2):281--308, 2008.

\bibitem{CS07}
C.~Galusinski and M.~Saad.
\newblock Two compressible immiscible fluids in porous media.
\newblock {\em J. Differential Equations}, 244:1741--1783, 2008.

\bibitem{ZS10}
Z.~Khalil and M.~Saad.
\newblock Solutions to a model for compressible immiscible two phase flow in
  porous media.
\newblock {\em Electronic Journal of Differential Equations}, 2010(122):1--33,
  2010.

\bibitem{ZS11}
Z.~Khalil and M.~Saad.
\newblock On a fully nonlinear degenerate parabolic system modeling immiscible
  gas-water displacement in porous media.
\newblock {\em Nonlinear Analysis}, 12:1591--1615, 2011.

\bibitem{kroener84}
D.~Kroener and S.~Luckhaus.
\newblock Flow of oil and water in a porous medium.
\newblock {\em J. Differential Equations}, 55:276--288, 1984.

\bibitem{kruzkov77}
S.~N. Kruzkov and S.~M. Sukorjanskii.
\newblock Boundary problems for systems of equations of two-phase porous flow
  type; statement of the problems, questions of solvability, justification of
  approximate methods.
\newblock {\em Math. USSR Sb.}, 33:62--80, 1977.

\bibitem{michel2003}
A.~Michel.
\newblock A finite volume scheme for the simulation of two-phase incompressible
  flow in porous media.
\newblock {\em SIAM J. Numer. Anal.}, 41:1301--1317, 2003.

\bibitem{peaceman}
D.W. Peaceman.
\newblock {\em Fundamentals of Numerical Reservoir Simulation}.
\newblock Elsevier Scientific Publishing, 1977.

\bibitem{these-bilal}
B.~Saad.
\newblock {\em Modélisation et simulation numérique d'écoulements
  multi-composants en milieu poreux}.
\newblock Thèse de doctorat de l'Ecole Centrale de Nantes, 2011.

\end{thebibliography}
\end{document}